\documentclass[matann2]{svjour}

\usepackage{amsmath}
\usepackage{amsfonts}
\usepackage{latexsym}
\usepackage{mathrsfs}
\usepackage{times}
\usepackage{array}
\usepackage{amscd}
\usepackage{a4}
\usepackage[mathcal]{euscript}
\usepackage{amssymb}
\usepackage{pstricks}
\usepackage{amsrefs}

\input{boxedeps} %% windows
\SetEPSFDirectory{} %% windows
\HideDisplacementBoxes
\SetRokickiEPSFSpecial  %% dvips by Tom Rokicki
\DeclareMathAlphabet{\ams}{U}{msb}{m}{n}

\def\sl{\text{SL}}

\def\gal{\text{Gal}}
\def\homeo{\text{Homeo}}
\def\id{\text{id}}

\def\rk{\text{rk}\,}
\def\coker{\text{coker}\,}
\def\im{\text{im}\,}

\def\ov{\overline}

\def\wtl{\widetilde}

\def\HH{\mathcal H}
\def\RR{\mathcal R}
\def\LL{\mathcal L}

\def\ve{\varepsilon}

\def\aa{\alpha}

\def\bb{\beta}
\def\ss{\sigma}

\def\ve{\varepsilon}

\def\Z{\ams{Z}}

\def\quo{/\kern -.45em\sim}

\newpsobject{showgrid}{psgrid}{subgriddiv=1,griddots=5,gridlabels=6pt}
\title{Galois theory, graphs and free groups}

\author{Brent Everitt\thanks{Some of the results of this paper were 
obtained while the author was visiting
the Department of Mathematics, University of Adelaide, Australia. 
}
}

\institute{\textsc{Brent Everitt}:
Department of Mathematics, University of York, York
YO10 5DD, United Kingdom. \email{bje1@york.ac.uk}
}

\titlerunning{}
\authorrunning{Brent Everitt 
}

\begin{document}

\maketitle

%\footnotetext{1991 {\em Mathematics Subject Classification:} 20F05, 30?} 

%%%%%%%%%%%%%%%%%%%%%%%%%%%%%%%%%%%%%%%%%%%%%%%%%%%%%%%%%%%%%%%%%%%%%%%%%

\begin{abstract}
A self-contained exposition is given of the topological and Galois-theoretic
properties of the category of combinatorial $1$-complexes, or graphs,
very much in the spirit of Stallings \cite{Stallings83}.
A number of classical, as well as some new results about free groups are 
derived.
\end{abstract}

%%%%%%%%%%%%%%%%%%%%%%%%%%%%%%%%%%%%%%%%%%%%%%%%%%%%%%%%%%%%%%%%%%%%%%%%%

\section*{Introduction}

This paper is about the interplay between graphs, free groups, and their
subgroups, a subject
with a long history that can be broadly divided into two schools. 
The first is combinatorial, where graphs, and particularly finite graphs,
provide a intuitively convenient way of picturing some aspects of the
theory of free groups,
as in for example \cites{Imrich77,Imrich76,Servatius83,Tardos96,Tardos92}.

The other approach is to treat graphs and their mappings 
as topological objects, a point of view with its origins from the very
beginnings of combinatorial
group theory, and resurrected in \cite{Stallings83} (see also \cites{Cohen89,
Gersten83,Neumann90}). This is the philosophy we take, but we differ from these 
earlier papers in that we place centerstage the theory of coverings of
arbitrary graphs, rather
than coverings being merely a prelude to immersions of finite graphs.
The first section sets up the topological preliminaries, 
\S \ref{topology:galoistheory} formulates the well known connection between
subgroups of free groups and coverings of graphs in a Galois-theoretic 
setting, while \S\S \ref{section:invariants}-\ref{section:pullbacks}
focus on the graph-theoretic implications
of finitely generated subgroups.

\section{The topology of graphs}\label{section:topological}

This section
is all very ``Stallings-esqe'' \cite{Stallings83}, with much of the material
in \S\S \ref{topological:graphs}-\ref{topological:pullbacksection}
well known. General references are \cites{Cohen89,Collins98,Gersten83,
Serre03,Scott79,Stallings83}. 
In \S \ref{topological:coverings} we deal with coverings, with a mixture
of well known and some (minor) new results; 
\S \ref{topology:lattice} introduces the lattice of intermediate coverings
of a cover.

\subsection{Graphs}\label{topological:graphs}

A {\em combinatorial $1$-complex\/} or {\em graph\/} \cite{Gersten83}*{\S 1.1} is 
an arbitrary set $\Gamma$ with an involutary map $^{-1}:\Gamma\rightarrow\Gamma$
and an idempotent map $s:\Gamma\rightarrow V_\Gamma$, (ie: $s^2=s$)
where $V_\Gamma$ is the
set of fixed points of $^{-1}$. Thus a graph has {\em vertices\/} $V_\Gamma$,
and {\em edges\/} $E_\Gamma:=\Gamma\setminus V_\Gamma$ with
(i). $s(v)=v$ for all $v\in V_\Gamma$;
(ii). $v^{-1}=v$ for all $v\in V_\Gamma$, $e^{-1}\in E_\Gamma$ and 
$e^{-1}\not= e=(e^{-1})^{-1}$ for all 
$e\in E_\Gamma$. Indeed, these two can be taken as a more transparent, but less 
elegant, definition. We will use both interchangebly.

The edge $e$ has {\em start vertex\/} $s(e)$ and {\em terminal vertex\/}
$t(e):=s(e^{-1})$; an {\em arc\/} is an edge/inverse
edge pair,
and an {\em orientation\/} for $\Gamma$ is a set $\mathcal{O}$ 
of edges containing exactly one edge from each arc. 
Write $\ov{e}$ for the arc
containing the edge $e$ (so that $\ov{e^{-1}}=\ov{e}$).
A pointed graph is a pair
$\Gamma_v:=(\Gamma,v)$ for $v\in\Gamma$ a vertex.

The graph $\Gamma$ is {\em finite\/} when $V_\Gamma$ is finite and 
{\em locally finite\/}
when the set $s^{-1}(v)$ is finite for every $v\in V_\Gamma$. 
The cardinality
of the set $s^{-1}(v)$ is the {\em valency\/} $\partial v$ of the vertex $v$.
A {\em path\/} is a finite sequence of edges, mutually incident in the obvious sense;
similarly we have {\em closed\/} paths and {\em trivial\/} paths (consisting of 
a single vertex).
$\Gamma$ is {\em connected\/} if any two vertices can be joined by a path. 
The connected component
of $\Gamma$ containing the vertex $v$ consists of those vertices for which 
there is a path connecting them to
$v$, together with all their incident edges.

\parshape=7 0pt\hsize 0pt.81\hsize 0pt.81\hsize 0pt.81\hsize 
0pt.81\hsize 0pt.81\hsize 0pt\hsize
A {\em map\/} of graphs is a set map $g :\Gamma\rightarrow \Lambda$
with $g(V_\Gamma)\subseteq V_\Lambda$, such that the diagram
\vadjust{\hfill\smash{\lower 62pt
\llap{
\begin{pspicture}(0,0)(2,2)
%\showgrid
\rput(-5.75,0.3){
\rput(0.7,-0.55){
\rput(5,2){$\Gamma$}\rput(6.55,2){$\Lambda$}
\rput(5,0.45){$\Gamma$}\rput(6.55,0.45){$\Lambda$}
\psline[linewidth=.1mm]{->}(5.3,2)(6.3,2)
\psline[linewidth=.1mm]{->}(5.3,0.45)(6.3,0.45)
\psline[linewidth=.1mm]{->}(5,1.7)(5,.7)
\psline[linewidth=.1mm]{->}(6.55,1.7)(6.55,.7)
\rput(4.7,1.2){$\ss_\Gamma$}\rput(6.85,1.2){$\ss_\Lambda$}
\rput(5.8,2.2){$g$}\rput(5.8,.65){$g$}
}}
\end{pspicture}}}}\ignorespaces
commutes, where $\ss_\Gamma$ is one of the $s$ or $^{-1}$ maps for $\Gamma$, 
and $\ss_\Lambda$ 
similarly, ie: $g s_\Gamma(x)=s_\Lambda g(x)$ and 
$g(x^{-1})=g(x)^{-1}$.
These are combinatorial versions of continuity: if $\Gamma$
is connected then $g(\Gamma)\subset\Lambda$ is connected.
A map is {\em dimension preserving\/} if $g(E_\Gamma)\subseteq E_\Lambda$.
These maps of graphs 
allow one to squash edges down to vertices as in \cite{Gersten83}, 
rather than the more rigid maps
of \cite{Serre03,Stallings83}. The pay off is that the quotient 
construction below is more useful.
A map $g:\Gamma_v\rightarrow\Lambda_u$ of pointed graphs is a map
$g:\Gamma\rightarrow\Lambda$ with $g(v)=u$.

A map $g:\Gamma\rightarrow \Lambda$ is a {\em homeomorphism\/} if it 
is dimension preserving and is a
bijection on the vertex and edge sets,
in which case
the inverse set map 
is a dimension preserving map of graphs
$g^{-1}:\Lambda\rightarrow \Gamma$, and hence a homeomorphism. 

The set of self
homeomorphisms $\Gamma\rightarrow \Gamma$
forms a group $\homeo(\Gamma)$, and a group action 
$G\stackrel{\varphi}{\rightarrow}\homeo(\Gamma)$
is said to {\em preserve
orientation\/} iff there is an orientation $\mathcal{O}$ for $\Gamma$ with 
$\varphi(g)(\mathcal{O})=\mathcal{O}$
for all $g\in G$. The action of $G$ is said to be 
{\em without inversions\/}
iff $\varphi(g)(e)\not= e^{-1}$ for all edges $e$ and for all 
$g\in G$. 
%
%\begin{lemma}
It is easy to see that $G$ 
preserves orientation if and only if it 
acts without inversions.
$G$ acts {\em freely\/} iff the action is free on the vertices, ie: 
for any $g\in G$ and $v$ a vertex, $\varphi(g)(v)=v$ implies that $g$
is the identity element. If $G$ acts freely and orientation preservingly, then
the action is free on the edges too.

A {\em subgraph\/} is a subset $\Lambda\subset\Gamma$,
such that the maps $s$ and $^{-1}$ give a graph when restricted to 
$\Lambda$. Equivalently, it is a graph mapping
$\Lambda\hookrightarrow\Gamma$ that is a homeomorphism
onto its image.
The {\em coboundary\/} $\delta\Lambda$ of a subgraph consists of those 
edges $e\in\Gamma$ with $s(e)\in\Lambda$ and $t(e)\not\in\Lambda$ 
(equivalently, it is those 
edges $e\in\Gamma$ with
$sq(e)$ the vertex $q(\Lambda)$
in the quotient complex $\Gamma/\Lambda$, where 
$q:\Gamma\rightarrow\Gamma/\Lambda$ is the quotient mapping as below).

An {\em elementary homotopy\/} of a path, 
$e_1\ldots e_ie_{i+1}\ldots e_k\leftrightarrow 
e_1\ldots e_i(ee^{-1})e_{i+1}\ldots e_k$
inserts or deletes a {\em spur\/}:
a path that consecutively traverses both edges of an arc $\overline{e}$.
Two paths %$w_1,w_2$ 
are (freely) {\em homotopic\/} %(written $w_1\simeq w_2$) 
iff there is a 
finite sequence of elementary homotopies
taking one to the other. Paths homotopic to a trivial path %based at 
are said to be {\em homotopically trivial\/}. It is easy to see that 
two homotopic paths have
the same start and terminal vertices (and thus homotopically 
trivial paths are
necessarily closed) and that homotopy is an equivalence relation on the
paths with common endpoints.

The {\em trivial graph\/} has a single vertex and no edges. 
The {\em real line\/} graph $\RR$ has vertices $V_\RR=\{v_k\}_{k\in\Z}$,
edges $E_\RR=\{e_k^{\pm 1}\}_{k\in\Z}$ and $s(e_k)=v_k,s(e_k^{-1})=v_{k+1}$.

\subsection{Quotients}

A {\em quotient relation\/} is an 
equivalence relation $\sim$ on $\Gamma$
such that
$$
\text{(i)}.\, x\sim y \Rightarrow s(x)\sim s(y)\text{ and }
x^{-1}\sim y^{-1},\hspace{2em}
\text{(ii)}.\, x\sim x^{-1}\Rightarrow[x]\cap V_\Gamma\not=\varnothing,
$$
where $[x]$ is the equivalence class under $\sim$ of $x$. 

\begin{proposition}%[quotient of a graph]
If $\sim$ is a quotient relation on a graph $\Gamma$ then the quotient graph 
$\Gamma\quo$ has vertices the equivalence classes $[v]$ for $v\in V_\Gamma$, 
edges the classes
$[e]$ for $e\in E_\Gamma$ with $[e]\cap V_\Gamma=\varnothing$, 
$[x]^{-1}=[x^{-1}]$, and
$s_{\Gamma/\sim}[x]=[s_\Gamma(x)]$.
Moreover, the quotient map
$q:\Gamma\rightarrow \Gamma\quo$ given by $q(x)=[x]$ is a map of graphs (and so
in particular, if $\Gamma$ is connected then $\Gamma\quo$ is connected).
\end{proposition}

Let $\Lambda_i\hookrightarrow\Gamma, (i\in I),$ be a set of mutually disjoint subgraphs
and define $\sim$ on $\Gamma$ by $x\sim y$ iff $x=y$ or both 
$x$ and $y$ lie in the same $\Lambda_i$. Write
$\Gamma/\Lambda_i$ for $\Gamma\quo$, the {\em quotient of $\Gamma$ by 
the family of subgraphs $\Lambda_i$}. It is what results by squashing
each $\Lambda_i$ to a distinct vertex. In particular, if the family consists
of a single subgraph $\Lambda\hookrightarrow\Gamma$, we have the quotient
$\Gamma/\Lambda$. The reader should be wary of the difference between the quotients
$\Gamma/\Lambda_i$ and $\Gamma/\Lambda$, for $\Lambda=\amalg\Lambda_i$ the union
of the disjoint subgraphs.

If $\sim$ is the equivalence
relation on $\Gamma$ consisting of the orbits of an action by a group $G$,
then $\sim$ is a 
quotient relation on $\Gamma$ if and only if the group action is orientation 
preserving. In this case we
may form the quotient complex $\Gamma/G:=\Gamma\quo$.

\subsection{Trees}\label{topological:trees}

A path in a graph is {\em reduced\/} when it contains no spurs;
by removing spurs,
any two vertices in the same component can be joined by a reduced path.

It is easily proved that 
for any vertices $u$ and $v$ of a graph $\Gamma$, there are $\leq 1$ reduced 
paths between them if and only if 
any closed non-trivial path in $\Gamma$ contains a spur
(equivalently, any closed path is homotopic to the trivial path based at 
one of its vertices).
A graph satisfying %either of the conditions of Proposition 
%\ref{uniqueness_reduced_paths_tree} 
any of these equivalent conditions is 
called a {\em forest\/}; a 
connected forest is a {\em tree\/}.

If $\Gamma$ is a finite graph with
$\partial v\geq 2$ for every vertex $v$, then it can be shown that
$\Gamma$ contains
a homotopically non-trivial closed path.
Hence if $T$ is a finite tree, then 
$|E_T|=2(|V_T|-1)$.

A {\em spanning forest\/} is a subgraph $\Phi\hookrightarrow\Gamma$
that is a forest
and contains all the vertices of $\Gamma$ (ie: $V_\Phi=V_\Gamma$). 
It is well known that 
%well ordering the edges of a connected $\Gamma$ that a 
spanning trees
can always be constructed for connected $\Gamma$.

\begin{proposition}\label{topological:trees:result100}
Let $T_i\hookrightarrow\Gamma$ be a family of mutually disjoint trees
in a connected graph $\Gamma$. Then there is a spanning tree $T\hookrightarrow
\Gamma$ containing the $T_i$ as subgraphs, and such that if $q:\Gamma\rightarrow
\Gamma/T_i$ is the quotient map, then $q(T)$ is a spanning tree for
$\Gamma/T_i$.
\end{proposition}

In particular, any spanning forest for $\Gamma$ can be extended to a spanning
tree.
For the proof, take $T$ to be $q^{-1}(T')$ for some spanning tree $T'$ of 
the (connected) graph $\Gamma/T_i$.

\subsection{The fundamental group}

The {\em fundamental group\/} $\pi_1(\Gamma,v)$ is the usual 
group of homotopy classes $[\gamma]$ of closed paths $\gamma$ at the 
vertex $v\in\Gamma$ (ie: equivalence classes under the homotopy relation) with
product $[\gamma_1][\gamma_2]=[\gamma_1\gamma_2]$. %, where 
If $\Phi$ is a forest, 
then $\pi_1(\Phi,v)$ is trivial for 
any vertex $v$ and conversely,
if $\Gamma$ connected has $\pi_1(\Gamma,v)$ trivial for some (hence every) 
vertex $v$, then $\Gamma$ is a tree.
A connected graph with trivial fundamental group is
{\em simply connected\/}.
A map $g:\Gamma_v\rightarrow \Lambda_u$ of graphs induces a group homomorphism
$g^*:\pi_1(\Gamma,v)\rightarrow\pi_1(\Lambda,u)$
by $g^*[\gamma]=[g(\gamma)]$ 
and this satisfies the usual functorality properties: 
$(\id)^*=\id$ 
and $(gf)^*=g^*f^*$.

\begin{proposition}%[excising trees]
\label{homotopy:excision}
If $T_i\hookrightarrow\Lambda$ is a family of mutually disjoint trees,
$v\in T\in\{T_i\}$ a vertex,
and $q:\Lambda\rightarrow\Lambda/T_i$ the quotient map, then
$q^*:\pi_1(\Lambda,v)\rightarrow 
\pi_1(\Lambda/T_i,q(v))$ is an isomorphism.
\end{proposition}

\begin{proof}
\parshape=8 0pt\hsize 0pt.8\hsize 0pt.8\hsize 0pt.8\hsize 0pt.8\hsize 
0pt.8\hsize 0pt.8\hsize 0pt\hsize
The key to the proof is that the quotient map $q$ is essentially just the 
identity map outside of 
the $T_i$.
\vadjust{\hfill\smash{\lower 20pt
\llap{\begin{pspicture}(0,0)(3,3)
\rput(-2.5,-2.3){
\rput(4,2){\BoxedEPSF{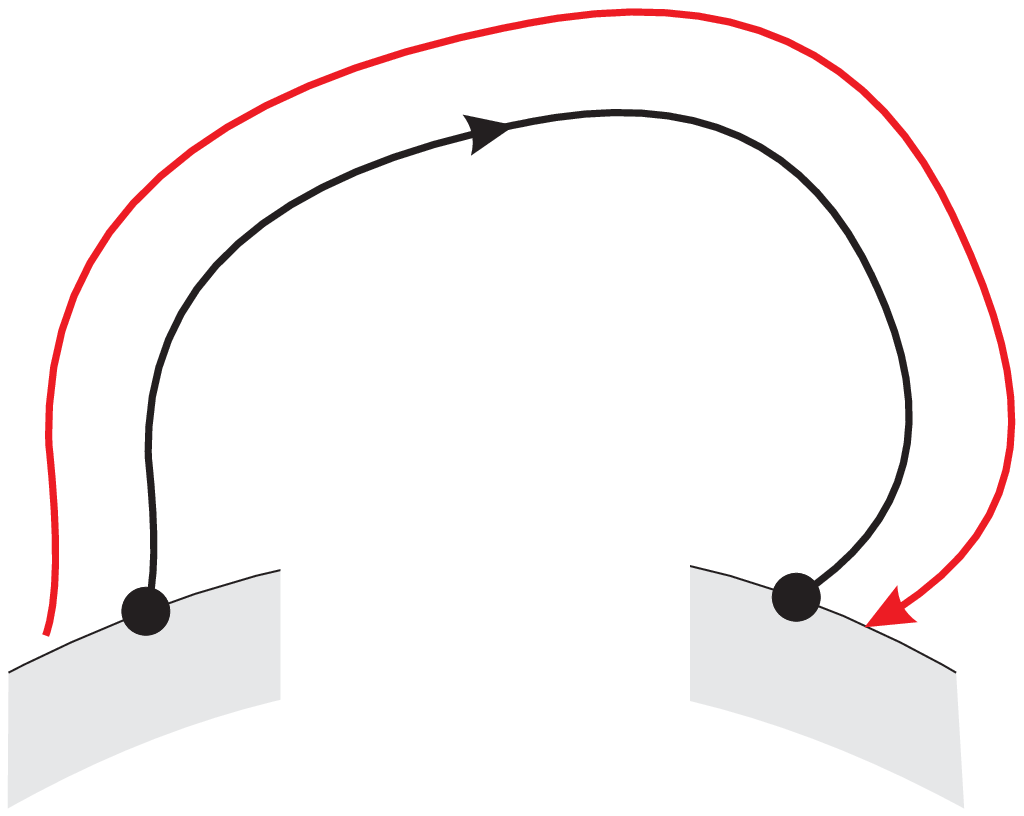 scaled 250}}
\rput(3.25,1.15){$T_{j-1}$}\rput(4.75,1.15){$T_{j}$}%\rput(4,3.2){$w'$}
%\rput(3.95,1.15){$v$}
\rput(5.35,2.6){{\red $\gamma_{j}'$}}
}
%\showgrid
\end{pspicture}}}}\ignorespaces
To see the surjectivity of $q^*$, suppose that $\gamma$ is a closed path in 
$\Lambda/T_i$ based at $q(v)$ 
and having edges $e_1\ldots e_k$. Then there are (unique) edges $e'_1,\ldots,
e'_k$ in $\Lambda$ with $q(e'_i)=e_i$ and 
$\gamma'=e'_1\ldots e'_k=\gamma'_1\gamma'_2\ldots\gamma'_k$ a sequence of paths
with the terminal vertex of $\gamma_{j}'$ and the start vertex of $\gamma'_{j+1}$
in the same tree $T_j$.
Use the connectedness of the $T_i$ to connect these up into 
a path in $\Lambda$ mapping via $q$ to $\gamma$.
For injectivity, suppose that $\gamma'$ is a closed path in $\Lambda$
based at $v$ and
mapping via $q$ to a homotopically trivial path $\gamma$ in $\Lambda/T_i$. 
If $\gamma$ contains a spur, then the section
of $\gamma'$ mapping to it looks %by Lemma \ref{images_of_paths}, 
like below:
$$
\begin{pspicture}(0,0)(14,1.5)
%\showgrid
\rput(0,0){
\rput(5,.8){\BoxedEPSF{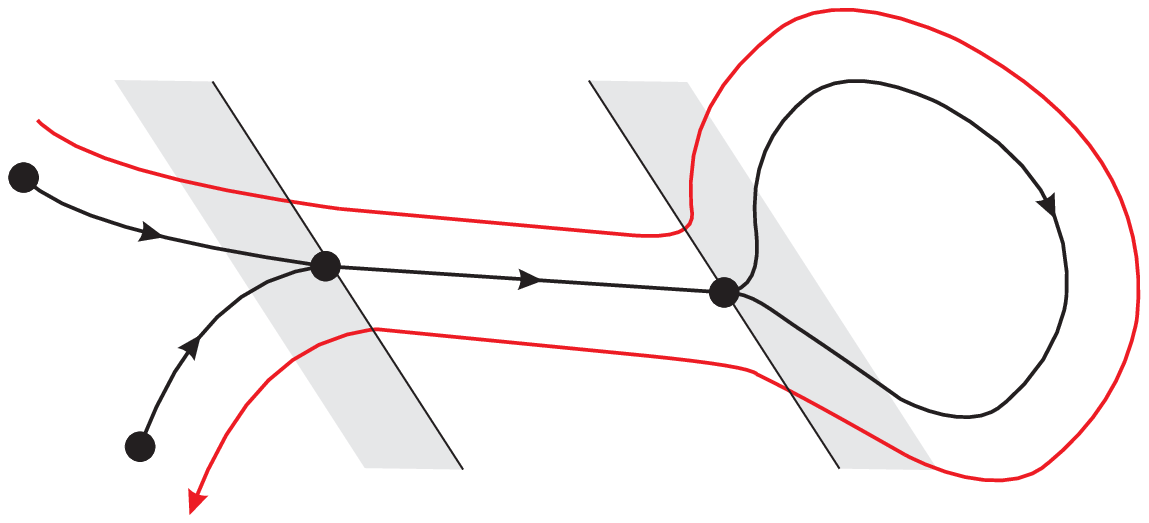 scaled 375}}
\rput(7.4,.5){{\red $\gamma'$}}
\rput(4.1,.05){$T_i$}\rput(5.3,1.5){$T_j$}
}
\psline[linewidth=.1mm]{->}(7.8,1)(9.2,1)
\rput(8.5,1.2){$q$}
\rput(0,0){
\rput(11,1){\BoxedEPSF{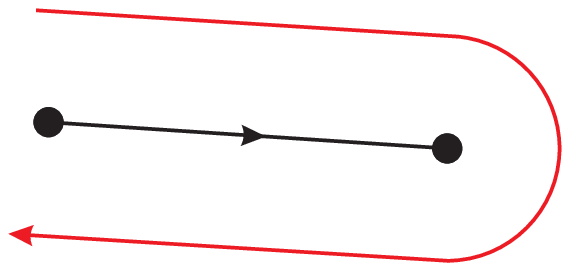 scaled 375}}
\rput(12.2,.5){{\red $\gamma$}}
}
\end{pspicture}
$$
Thus a sequence of elementary homotopies reducing $\gamma$ to the trivial 
path in $\Lambda/T_i$ can be mirrored by
homotopies in $\Lambda$ that reduce $\gamma'$ to a closed path completely 
contained in $T$. As $T$
is simply connected, this path can in turn be homotoped to the trivial path. 
Thus, only homotopically 
trivial paths can be sent by $q$ to homotopically trivial paths, so $q^*$ is 
injective. 
\qed
\end{proof}

Fix a spanning tree $T\hookrightarrow\Gamma$, choose an edge $e$ from each arc of
$\Gamma$, 
and consider the homotopy class of the path through $T$ from 
$v$ to $s(e)$, traverses $e$ and travels back through $T$ to $v$. Then 
{\em Schreier generators\/} for
$\pi_1(\Gamma,v)$ are the homotopy classes of such paths arising from 
the arcs omitted by $T$.

\subsection{Homology}

Fix an orientation $\mathcal{O}$ for $\Gamma$,
and always write arcs in the form $\ov{e}$ for 
$e\in \mathcal{O}$, and paths in the form 
$\gamma=e_1^{\ve_1}\ldots e_k^{\ve_k}$ with
$e_i\in \mathcal{O}$ and $\ve_i=\pm 1$. Let $\Z[V_\Gamma]$ and 
$\Z[\text{arcs}]$ be the free abelian groups on the
vertices and arcs of $\Gamma$ (alternatively, one can take
$\Z[E_\Gamma]$
and then pass to the quotient $\Z[E_\Gamma]/\langle e+e^{-1}=0\rangle$; 
we prefer the more concrete version).
Define the boundary of an arc $\ov{e}$ to be 
$\partial(\ov{e})=t(e)-s(e)\in\Z[V_\Gamma]$, and for 
$\sum n_i\ov{e}_i\in\Z[\text{arcs}]$, let $\partial(\sum
n_i\ov{e}_i)=\sum n_i\partial(\ov{e}_i)$.
Then $\partial$ is a group homomorphism 
$\partial:\Z[\text{arcs}]\rightarrow\Z[V_\Gamma]$, and
the {\em homology\/} of $\Gamma$ is the pair of abelian
groups
$$
H_1(\Gamma)=\ker\partial\text{ and }H_0(\Gamma)=\coker\partial,
$$
(ie: $H_0(\Gamma)=\Z[V_\Gamma]/\im\partial$).

By following the proofs in the topological category,
%(eg: \cite{Hatcher02}), 
one can show the standard homological facts:
$H_0(\Gamma)$ is free abelian on the connected components of $\Gamma$;
if $\Gamma$ is single vertexed, then $H_1$ is free abelian on the arcs.
If $\gamma=e_1^{\ve_1}e_2^{\ve_2}\ldots e_k^{\ve_k}$ is a closed path
at $v$ then $\partial(\sum \ve_i\ov{e}_i)=0$, and the Hurewicz map sending
$\gamma$ to $\sum \ve_i\ov{e}_i$ is well defined upto homotopy,
thus, for $\Gamma$ connected,
a surjective homomorphism $\pi_1(\Gamma,v)\rightarrow H_1(\Gamma)$
with kernel the commutator 
subgroup of $\pi_1(\Gamma,v)$.

In particular, $H_1(\Gamma)$ is the abelianisation of $\pi_1(\Gamma,v)$, 
so that if $\Gamma_1,\Gamma_2$ are connected graphs
with $\pi_1(\Gamma_1,v_1)\cong\pi_1(\Gamma_2,v_2)$ then 
$H_1(\Gamma_1)\cong H_1(\Gamma_2)$.

\subsection{Rank and spines}\label{topological:rank}

Graph homology provides an important invariant for graphs:

\begin{proposition}[rank an invariant]
\label{topological:rank:result100}
Let $T\hookrightarrow\Gamma$ be a spanning tree for $\Gamma$ connected.
Then $H_1(\Gamma)$ is free abelian 
with basis the set of arcs of $\Gamma$ omitted by $T$. 
\end{proposition}

Thus the cardinality of the set of omitted arcs is equal to 
$\rk_\Z H_1(\Gamma)$ and independent of $T$
(this can also be
shown directly without recourse to homology). 
Define the
{\em rank\/} $\rk\Gamma$ of $\Gamma$ connected to be $\rk_\Z H_1(\Gamma)$, 
or the cardinality of the set of arcs omitted by a spanning tree. 
%For an arbitrary graph, define the rank to be the cardinality of the set
%of arcs omitted by spanning trees for each connected component.

\begin{proof}[of Proposition \ref{topological:rank:result100}]
We have $\pi(\Gamma,v)\cong\pi_1(\Gamma/T,q(v))$ by Proposition 
\ref{homotopy:excision}, 
hence $H_1(\Gamma)\cong H_1(\Gamma/T)$,
%by Proposition \ref{homology100} part (4), 
with $\Gamma/T$ single vertexed, hence this final homology free 
abelian on its arcs, 
%by Proposition \ref{homology100} part (2), 
ie: free abelian on the
arcs of $\Gamma$ omitted by $T$.
\qed
\end{proof}

If $\Gamma$ is finite, locally finite, connected, then,
$2(\rk\Gamma-1)=
|E_\Gamma|-2|V_\Gamma|$ by
\S \ref{topological:trees}; 
clearly, $\rk\Gamma=0$ if and only if $\Gamma$ is a tree.
If $\Gamma$ a connected graph and $T_i\hookrightarrow\Gamma$ a set of 
mutually disjoint trees, then $\rk\Gamma=\rk\Gamma/T_i$ (this follows either
from Proposition \ref{homotopy:excision} using $\rk=\rk H_1$ or by Proposition 
\ref{topological:trees:result100} using rank 
the number of arcs omitted by a spanning tree).

If $\Lambda$ is a connected graph and $v$ a vertex, 
then the {\em spine $\widehat{\Lambda}_v$ of $\Lambda$ at $v$\/}, 
is defined
to be the union in $\Lambda$ of all closed reduced paths starting at $v$.
Stallings and others use core graphs; we have followed \cite{Neumann90}.

\begin{lemma}\label{topological:spines:result100}
%\begin{enumerate}
(i). $\widehat{\Lambda}_v$ is connected with 
$\rk\widehat{\Lambda}_v=\rk\Lambda$.
(ii). If $u\in\widehat{\Lambda}_v$, then every closed reduced path
starting at $u$ is contained in $\widehat{\Lambda}_v$.
(iii). Spines are topological invariants, ie: 
a homeomorphism $f:\Lambda_u\rightarrow\Delta_v$
restricts to a homeomorphism $\widehat{\Lambda}_u\rightarrow\widehat{\Delta}_v$.
\end{lemma}

\begin{proof}
The connectedness is immediate. 
If $T$ is a spanning tree for $\Lambda$ and $e$ an edge not in $T$, then
$e$ is contained in the spine $\widehat{\Lambda}_v$, for,
the closed path obtained by traversing the reduced path
through $T$ from $v$
to $s(e)$, across $e$ and back via the reduced path through $T$
is reduced. The rank assertion follows.
For part (ii), let $\mu$ be a closed reduced path at $u$ and $\gamma$ a 
reduced path in the spine from $v$ to $u$. Then $\gamma=\gamma'_1\gamma_1=
\gamma'_2\gamma_2$ and $\mu=\gamma_1^{-1}\mu'\gamma_2$ with 
$\gamma'_1\mu'(\gamma'_2)^{-1}$ reduced at $v$, hence in the spine.
A homeomorphism sends closed reduced paths to closed reduced paths
(compare with Proposition \ref{topological:coverings:result100} and the maps of 
\S \ref{topological:coverings}), 
hence 
$f(\widehat{\Lambda}_v)\subset\widehat{\Delta}_u$, and the converse similarly
using $f^{-1}$.
\qed
\end{proof}

\subsection{Pushouts}\label{topological:pushouts}

These are important examples of quotients. Let
$\Lambda_1,\Lambda_2$ and $\Delta$ be graphs and 
$g_i:\Delta\rightarrow \Lambda_i$ 
maps of graphs. 
Let $\sim$ on the disjoint union $\Lambda_1\amalg\Lambda_2$ be the equivalence
relation {\em generated by\/} the $x\sim y$ iff there is a $z\in\Delta$ with
$x=g_1(z)$ and $y=g_2(z)$. Thus, $x\sim y$ iff
there are $x_0,x_1,\ldots,x_k$ with $x_0=x$ and $x_k=y$, and for each $j$, 
there is $z\in \Delta$ with $x_j=g_i(z)$, $x_{j+1}=g_{i+1\text{ mod }2}(z)$.
If $\sim$ is a 
quotient
relation then call the quotient $\Lambda_1\coprod \Lambda_2\quo$ the 
{\em pushout\/}
of the data $g_i:\Delta\rightarrow \Lambda_i$,
denoted 
$\Lambda_1\coprod_\Delta \Lambda_2$.

Given graphs and maps as above, 
the pushout cannot always
be formed, precisely because the quotient cannot always be formed. 
%
%\begin{lemma}[\cite{Stallings83}]
Stallings \cite{Stallings83}*{page 552} shows that if the $g_i$ are dimension
preserving,
then the pushout exists if and only if there are orientations
$\mathcal{O},\mathcal{O}_i$ for $\Delta,\Lambda_i$
with $g_i(\mathcal{O})\subseteq \mathcal{O}_i$.
Thus in particular, if the graphs $g_1(\Delta)$ and $g_2(\Delta)$ 
are disjoint, then the pushout can always be formed.

Define $t_i:\Lambda_i\rightarrow\Lambda_1\coprod_\Delta\Lambda_2$ 
to be the compositions $\Lambda_i\hookrightarrow 
\Lambda_1\coprod \Lambda_2
\rightarrow \Lambda_1\coprod \Lambda_2\quo$ of the 
inclusion of $\Lambda_i$ in the disjoint union and the quotient map.

\begin{proposition}\label{pushout:categorical}
If $\Delta\not=\varnothing$ and the $\Lambda_i$ are connected then the pushout
is connected, and the maps $t_i$ make the
diagram on the left commute. 
$$
\begin{pspicture}(6,3)
%\showgrid
\rput(0,-.5){
\rput(-1.5,1){
\rput(0,2){$\Delta$}
\rput(0,0){$\Lambda_1$}
\rput(2,2){$\Lambda_2$}
\rput(2.1,0){$\Lambda_1\coprod_\Delta \Lambda_2$}
\psline[linewidth=.1mm]{->}(0,1.7)(0,.3)
\psline[linewidth=.1mm]{->}(.3,2)(1.7,2)
\psline[linewidth=.1mm]{->}(.3,0)(1.2,0)
\psline[linewidth=.1mm]{->}(2,1.7)(2,.3)
\rput(.25,1.05){$g_1$}
\rput(1,1.8){$g_2$}
\rput(.75,.2){$t_1$}
\rput(1.8,1){$t_2$}
}
\rput(0,.5){
\rput(5,1){
\rput(0,2){$\Delta$}
\rput(0,0){$\Lambda_1$}
\rput(2,2){$\Lambda_2$}
\rput(2.1,0){$\Lambda_1\coprod_\Delta \Lambda_2$}
\psline[linewidth=.1mm]{->}(0,1.7)(0,.3)
\psline[linewidth=.1mm]{->}(.3,2)(1.7,2)
\psline[linewidth=.1mm]{->}(.3,0)(1.2,0)
\psline[linewidth=.1mm]{->}(2,1.7)(2,.3)
\rput(.25,1.05){$g_1$}
\rput(1,1.8){$g_2$}
\rput(.75,.2){$t_1$}
\rput(1.8,1){$t_2$}
}
\rput(8,0){$B$}
\psbezier[linewidth=.1mm]{->}(5.2,.8)(6,0)(7,0)(7.8,0)
\psbezier[linewidth=.1mm]{->}(7.2,2.8)(8,2)(8,1)(8,.2)
\psline[linewidth=.1mm]{->}(7.3,.7)(7.8,.2)
\rput(6.35,.4){$t'_1$}\rput(8,2){$t'_2$}
}}
\end{pspicture}
$$
Moreover the pushout is universal in that if $B$,
$t'_1$, $t'_2$ are a graph and maps making 
such a square commute, then there
is a map $\Lambda_1\coprod_\Delta \Lambda_2\rightarrow B$ making the 
diagram above right commute.
\end{proposition}

These properties can of course be taken as an alternative, categorical definition
of the pushout, with uniqueness following from the universality and the usual
formal nonsense.
If the $g_i:\Delta_v\rightarrow\Lambda_{u_i}$ are pointed maps, 
and $q:\Lambda_1\coprod \Lambda_2\rightarrow \Lambda_1\coprod \Lambda_2\quo$
the quotient map, then $q(u_1)=q(u_2)=u$ (say), and we have a pointed
version of Proposition \ref{pushout:categorical}, involving the pointed pushout
$(\Lambda_1\coprod_\Delta\Lambda_2)_u$.

Many of the quotient constructions from topology (eg: cone, 
suspension, $\ldots$) 
can be expressed as some kind of pushout or other, 
%, and there are other important 
%examples 
%(eg: Stallings fold). 
but we content
ourselves with the following:
let $\Delta$ be a graph with $E_\Delta=\varnothing$ and the 
$g_i:\Delta\hookrightarrow\Lambda_i$ homeomorphisms
onto their images (ie: injections on the vertices of $\Delta$). 
The resulting pushout (which always exists),
the {\em wedge sum\/} $\Lambda_1\bigvee_{\Delta}\Lambda_2$,
is the result of identifying the vertices
of copies of $\Delta$ in the $\Lambda_i$. If the $\Lambda_i$ coincide
($=\Lambda$ say) with maps $\Delta\rightrightarrows\Lambda$, 
then we write $\bigvee_\Delta\Lambda$.

If $\Delta$ is the trivial graph, $T$ a tree and 
$g_i:\Delta\rightrightarrows T$ distinct maps,
then the wedge sum $\bigvee_{\Delta}T$ has a non-trivial reduced 
closed path that is unique upto cyclic reordering. Thus,
by removing a single arc from $\bigvee_{\Delta}T$ we obtain a new tree.

\subsection{Pullbacks}\label{topological:pullbacksection}

The categorical nature of the pushout construction 
(ie: Proposition \ref{pushout:categorical}) 
suggests a ``co-'' version:
let $\Lambda_1,\Lambda_2$ and $\Delta$ be graphs and 
$g_i:\Lambda_i\rightarrow \Delta$ 
maps of graphs.
The {\em pullback\/} $\Lambda_1\prod_\Delta \Lambda_2$  has vertices
(resp. edges) the $x_1\times x_2$, $x_i\in V_{\Lambda_i}$ (resp.
$x_i\in E_{\Lambda_i}$) such that $g_1(x_1)=g_2(x_2)$, and
$s(x_1\times x_2)=s(x_1)\times s(x_2)$, $(x_1\times x_2)^{-1}
=x_1^{-1}\times x_2^{-1}$. 
Taking $\Delta$ to be the trivial graph has the 
effect of removing the 
$g_1(x)=g_2(y)$
conditions and 
the result is the {\em product\/} 
$\Lambda_1\prod \Lambda_2$. Thus
the pullback $\Lambda_1\prod_\Delta \Lambda_2$ is a subgraph of the 
product $\Lambda_1\prod\Lambda_2$,
but the product will have many more vertices and edges.
Define maps $t_i:\Lambda_1\prod_\Delta\Lambda_2\rightarrow\Lambda_i$ to be
the compositions $\Lambda_1\prod_\Delta\Lambda_2\hookrightarrow
\Lambda_1\prod\Lambda_2\rightarrow\Lambda_i$, with the second map the projection
$x_1\times x_2\mapsto x_i$.

\begin{proposition}\label{topological:pullbacks}
The $t_i$ are 
dimension preserving maps
making the diagram below left
commute,
$$
\begin{pspicture}(6,3)
%\showgrid
\rput(0,-1){
\rput(-1.5,1.5){
\rput(0,2){$\Lambda_1\prod_\Delta \Lambda_2$}
\rput(0,0){$\Lambda_1$}
\rput(2,2){$\Lambda_2$}
\rput(2,0){$\Delta$}
\psline[linewidth=.1mm]{->}(0,1.7)(0,.3)
\psline[linewidth=.1mm]{->}(.9,2)(1.7,2)
\psline[linewidth=.1mm]{->}(.3,0)(1.7,0)
\psline[linewidth=.1mm]{->}(2,1.7)(2,.3)
\rput(.25,1){$t_1$}
\rput(1.3,1.8){$t_2$}
\rput(1,.2){$g_1$}
\rput(1.75,1){$g_2$}
}
\rput(5,1){
\rput(0,2){$\Lambda_1\prod_\Delta \Lambda_2$}
\rput(0,0){$\Lambda_1$}
\rput(2,2){$\Lambda_2$}
\rput(2,0){$\Delta$}
\psline[linewidth=.1mm]{->}(0,1.7)(0,.3)
\psline[linewidth=.1mm]{->}(.9,2)(1.7,2)
\psline[linewidth=.1mm]{->}(.3,0)(1.7,0)
\psline[linewidth=.1mm]{->}(2,1.7)(2,.3)
\rput(.25,1){$t_1$}
\rput(1.3,1.8){$t_2$}
\rput(1,.2){$g_1$}
\rput(1.75,1){$g_2$}
}
\rput(4,4){$B$}
\psbezier[linewidth=.1mm]{->}(4,3.8)(4,3)(4,2)(4.8,1.2)
\psbezier[linewidth=.1mm]{->}(4.2,4)(5,4)(6,4)(6.8,3.2)
\psline[linewidth=.1mm]{->}(4.2,3.8)(4.7,3.3)
\rput(3.8,2.6){$t'_1$}\rput(5.6,4.1){$t'_2$}
}
\end{pspicture}
$$
Moreover, the pullback is universal in that if $B$,
$t'_1,t'_2$ are a graph and maps making such a square commute,
then
there is a map $B\rightarrow \Lambda_1\prod_\Delta \Lambda_2$ making the
diagram above right commute.
\end{proposition}

In general the pullback need not be connected.
If the $g_i:\Lambda_{u_i}\rightarrow\Delta_v$ are pointed maps then 
$u_1\times u_2$ is a vertex of the pullback, and we may consider the
connected component containing $u_1\times u_2$. Call this the 
{\em pointed pullback\/} $(\Lambda_1\prod_\Delta\Lambda_2)_{u_1\times u_2}$,
and we then have a pointed version of Proposition 
\ref{topological:pullbacks}.
In most of our usages of the pullback construction, the graph $\Delta$ 
will be single vertexed, and so the 
vertex set will just be $V_{\Lambda_1}\times V_{\Lambda_2}$.

\subsection{Coverings}\label{topological:coverings}

A map %$\widetilde{K}\stackrel{p}{\longrightarrow} K$
%\marginlabel{\small{$\blob$ pointed coverings}}
$p:\Lambda\rightarrow\Delta$ of graphs is a {\em covering\/} iff
(i). $p$ preserves dimension; and
%(ie: maps vertices to vertices, edges to edges and faces to faces);
(ii). for every vertex $v\in \Lambda$, $p$ is a bijection from the set of 
edges in $\Lambda$
with start vertex $v$ to the set of edges in $\Delta$ with start
vertex $p(v)$.
If
$p(x)=y$, then one says that $x$ {\em covers\/} $y$,
and $y$ {\em lifts\/} to $x$.
The set of all lifts of the cell $y$, or the set of all cells covering
$y$, is its {\em fiber\/} $\text{fib}_{\Lambda\rightarrow\Delta}(y)$. 

\begin{proposition}[\cite{Stallings83}*{\S 4.1}]
\label{topological:coverings:result100}
%\marginlabel{\small{$\blob$ ie: $p:(\Lambda,u)\rightarrow (\Delta,v)$ a
%covering}}
Let $p:\Lambda_u\rightarrow\Delta_v$ be a covering. 
%and $u$ a vertex of $\Lambda$ 
%with $p(u)=v$.
%\begin{enumerate}

\noindent(i). If $\gamma$ is a path in $\Delta$ starting at $v$
then there is a path $\gamma'$ in $\Lambda$ starting at $u$
and covering $\gamma$.
%\noindent(ii). 
Moreover, if $\gamma_1,\gamma_2$ are paths in $\Lambda$ starting at 
$u$ and 
%$p(\wtl{\gamma}_1)=p(\wtl{\gamma}_2)$, 
covering the same path, 
then $\gamma_1=\gamma_2$.

\noindent(ii). A path in $\Lambda$ covering a spur is itself a spur.
Consequently, two paths in $\Lambda$ 
covering homotopic paths are homotopic.

\noindent(iii). If $g:\Gamma_w\rightarrow\Delta_v$ is a map then there
is a map $f:\Gamma_w\rightarrow\Lambda_u$ with $g=fp$ if and only
if $g^*\pi_1(\Gamma,w)\subset p^*\pi_1(\Lambda,u)$.

\noindent(iv). $p^*:\pi_1(\Lambda,u)\rightarrow\pi_1(\Delta,v)$ is injective, and if
$u'$ is the terminal vertex of a path $\mu$ starting at $u$, then 
$p^*\pi_1(\Lambda,u)=gp^*\pi_1(\Lambda,u')g^{-1}$, where $g$ is the 
homotopy class of $p(\mu)$.
%\end{enumerate}
\end{proposition}

The path $\gamma'$ in (i) is a lift of $\gamma$ to $u$, such lifts being
unique by (ii). The combination of these two is called {\em path lifting\/}, 
while (ii) is  {\em spur-lifting\/} and {\em homotopy lifting\/}. 
In particular, the image under a covering of a reduced path is reduced (whereas,
as spurs always map to spurs, the pre-image of a reduced path is 
reduced under any mapping). Part (iii) is a general lifting criterion that implies
in particular that
if $\gamma$ a closed path at $v$ then its homotopy class lies in
$p^*\pi_1(\Lambda,u)$ if and only if there is a closed path $\mu$ at $u$
with $p(\mu)=\gamma$. Part (iv) follows immediately from this and homotopy lifting.

\begin{lemma}
\label{topological:coverings:result200}
Let $p:\Lambda\rightarrow\Delta$ be a covering. 
%\begin{enumerate}

\noindent(i). If $\Delta$ is connected then $p$ maps the cells of $\Lambda$
surjectively onto the cells of $\Delta$.

\noindent(ii). If $\Lambda$ is connected then 
the fibers of any two cells of $\Delta$ have the same cardinality, called the degree, 
$\deg(\Lambda\rightarrow\Delta)$, of the covering.

\noindent(iii). If $\Lambda,\Delta$ are connected and 
$\deg(\Lambda\rightarrow\Delta)=1$, 
then the covering 
$\Lambda\rightarrow\Delta$ is a
homeomorphism.
%\end{enumerate}
\end{lemma}

\begin{proof}
In (i), surjectivity on the vertices follows by path lifting and on the
edges by definition. 
Path lifting gives a bijection in (ii) between the fibers of two vertices,
and between the fiber of an edge and it's start vertex.
Part (iii) follows immediately from (i) and (ii).
\qed
\end{proof}

From now on, {\em all\/} coverings will be maps between connected complexes
unless stated otherwise.

\begin{lemma}\label{topology:coverings:result300}

\noindent(i). Let
$\Lambda\stackrel{q}{\rightarrow}\Gamma\stackrel{r}{\rightarrow}\Delta$ be
maps with $p=rq$. 
If any two of $p,q$ and $r$ are coverings, then so is the third.

\noindent(ii). If a group $G$ acts orientation preservingly and freely on $\Lambda$
then the quotient map $q:\Lambda\rightarrow\Lambda/G$ is a covering.
\end{lemma}

Call the coverings $\Lambda\rightarrow\Gamma\rightarrow\Delta$ in (i) 
{\em intermediate\/} to the covering $\Lambda\rightarrow\Delta$. It follows from the
comments following Proposition \ref{topological:coverings:result100} that
if $p:\Lambda_u\rightarrow\Delta_v$ and $r:\Gamma_x\rightarrow\Delta_v$, then
 $p^*\pi_1(\Lambda,u)\subset r^*\pi_1(\Gamma,x)$.

\begin{proof}[of Lemma \ref{topology:coverings:result300}]
The freeness of the action in (ii) ensures the injectivity of $q$ on the edges
starting at a vertex of $\Lambda$. Part (i) is an easy exercise.
\qed
\end{proof}

\begin{proposition}\label{topology:coverings:result400}
%\begin{enumerate}
Let $\Lambda$ be a graph and $\Upsilon_1,\Upsilon_2\hookrightarrow\Lambda$
subgraphs of the form,
$$
\begin{pspicture}(14,1)
%\showgrid
\rput(4,0){
\rput(3,.5){\BoxedEPSF{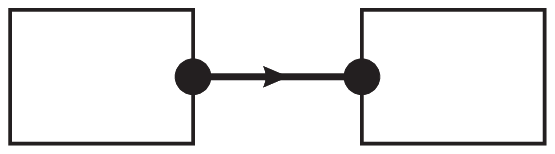 scaled 750}}
\rput(.5,.5){$\Lambda=$}%\rput(2.5,.4){$u$}
\rput(1.6,.5){$\Upsilon_1$}\rput(4.4,.5){$\Upsilon_2$}
\rput(3,.7){$e$}
}
\rput(14,0.5){$(\dag)$}
\end{pspicture}
$$
\noindent(i). If $p:\Lambda\rightarrow\Delta$ is a covering with $\Delta$ 
single vertexed, then the real line is a subgraph 
$\aa:\RR\hookrightarrow\Lambda$,
with $\aa(e_0)=e$ and $p\aa(e_k)=p(e)$ for all $k\in\Z$.

\noindent(ii). If $\Upsilon_1$ is a tree, $p:\Lambda\rightarrow\Delta$,
$r:\Gamma\rightarrow\Delta$ coverings, and
$\aa:\Upsilon_2\hookrightarrow\Gamma$ a homeomorphism onto its image,
then
there is an intermediate
covering $\Lambda\stackrel{q}{\rightarrow}\Gamma\stackrel{r}{\rightarrow}\Delta$.

\noindent(iii). If $\Psi\rightarrow\Lambda$ is a covering and $\Upsilon_1$ a tree,
then $\Psi$ also has the form $(\dag)$ for some subgraphs 
$\Upsilon'_1,\Upsilon'_2\hookrightarrow\Psi$, with $\Upsilon'_1$ a tree.
%\end{enumerate}
\end{proposition}

\begin{proof}
(i). Lift the edge $p(e)$ to the vertex $t(e)$ to get an edge $e_1$ of $\Lambda$. 
The form of $\Lambda$ 
prohibits $t(e_1)$ from being any vertex of $\Upsilon_1$, except possibly 
$s(e)$, in which
case $e_1=e^{-1}$. But then
$p(e)^{-1}=p(e^{-1})=p(e_1)=p(e)$ a contradiction. Thus $e_1$ is an edge 
and $t(e_1)$ a vertex
of $\Upsilon_2$, and if $t(e_1)=t(e)$, then the injectivity of $p$ fails 
at this common vertex
(as then both $e^{-1}$ and $e_1^{-1}$ start at $t(e_1)$ and cover 
$p(e)^{-1}$). We therefore have $t(e_1)\not=t(e)$ and this process can be
continued inductively, giving the ``positive'' half of $\RR$ a subgraph of 
$\Upsilon_2$. The symmetry of $\Lambda$ gives the negative half a subgraph
of $\Upsilon_1$.

For (ii), it suffices, by part Lemma \ref{topology:coverings:result300}(i), 
to find a map $q:\Lambda\rightarrow\Gamma$
with $p=rq$. Let $q$ coincide with $\aa$ on $\Upsilon_2$. For any vertex of 
$\Upsilon_1$,
take the reduced path to it from $t(e)$, project via $p$ to $\Delta$, and lift to
$\aa t(e)\in\Gamma$. The edges of $\Upsilon_1$ (and $e$) are similar.

(iii). Let $v=t(e)$ and $u$ be in the fiber of $v$ via the covering 
$\Psi\rightarrow\Lambda$. Take a reduced path in
$\Upsilon_1\cup\{e\}$ from $v$ to each vertex of this tree and lift to a path
at $u$.
Let $\Upsilon'_1$ be the union in $\Psi$
of these lifted paths. A closed path in $\Upsilon'_1$ at $u$  covers
a closed path at $v$ in $\Upsilon_1$, a tree, hence by spur-lifting,
$\Upsilon'_1$ is a tree. If $e_1,e_2$ are edges in the 
coboundary $\delta\Upsilon'_1$ then they cover edges in the coboundary
$\delta\Upsilon_1$, ie: they cover $e$. A reduced path in $\Upsilon'_1$
from $s(e_1)$ to $s(e_2)$ covers a reduced closed path in
$\Upsilon_1$ at $s(e)$. As this covered path must be trivial we get
$s(e_1)=s(e_2)$, hence $e_1=e_2$. Thus $\Upsilon'_1$ has a single coboundary
edge as required.
\qed
\end{proof}

\begin{proposition}\label{topology:coverings:result500}
Let $p:\Lambda\rightarrow\Delta$ be a covering %with $\Delta$ connected 
and $T\hookrightarrow\Delta$ a 
tree. Then
(i). $p^{-1}(T)$ a forest. (ii). If $T_i\hookrightarrow\Lambda, (i\in I)$ are the 
connected components of $p^{-1}(T)$, then $p$ maps each $T_i$
homeomorphically onto $T$. (iii). There is an induced
covering $\Lambda/T_i\rightarrow\Delta/T$ making the
diagram,
$$
\begin{pspicture}(2,2)
%\showgrid
\rput(0,0){
\rput(0,2){$\Lambda$}
\rput(0,0){$\Delta$}
\rput(2,2){$\Lambda/T_i$}
\rput(2,0){$\Delta/T$}
\psline[linewidth=.1mm]{->}(0,1.7)(0,.3)
\psline[linewidth=.1mm]{->}(.3,2)(1.55,2)
\psline[linewidth=.1mm]{->}(.3,0)(1.55,0)
\psline[linewidth=.1mm]{->}(2,1.7)(2,.3)
%\rput(-.175,1.05){$p$}
%\rput(1,1.8){$q$}
%\rput(1,.2){$q$}
%\rput(2.2,1.075){$p'$}
}
\end{pspicture}
$$
commute (where the horizontal maps are the quotients) and such that
$\deg(\Lambda/T_i\rightarrow\Delta/T)=\deg(\Lambda\rightarrow\Delta)$.
\end{proposition}

This procedure is independent of the tree $T$: if $T'\hookrightarrow\Delta$
another tree such that there is a homeomorphism $\aa:\Delta/T\rightarrow\Delta/T'$
with $\aa q=q'$ for $q,q':\Delta\rightarrow\Delta/T,\Delta/T'$ the quotient
maps, then by Proposition \ref{topological:coverings:result100}(iv), 
there is a homeomorphism $\Lambda/T_i\rightarrow\Lambda/T'_i$. Typically
we will take $T,T'$ to be spanning trees, so that $\Delta/T,\Delta/T'$ are single
vertexed with $\rk\Delta$ loops, and such an $\aa$ is easily found.

\begin{proof}[of Proposition \ref{topology:coverings:result500}]
That $p^{-1}(T)$ is a forest follows by spur-lifting.
For (ii), $p$ is injective on the vertices (and hence edges) of $T_i$ as $T$ is a tree
and by spur lifting;
surjectivity follows by path lifting.
If $q',q$ are the top and bottom quotient maps, define
for any cell $q'(x)\in\Lambda/T_i$, the map $p'q'(x)
=qp(x)$. Taking a vertex $v\in T$, the degree assertion follows immediately from
(ii).
\qed
\end{proof}

A covering $p:\wtl{\Delta}_u\rightarrow\Delta_v$ is {\em universal\/} iff
for any covering $r:\Gamma_w\rightarrow\Delta_v$ there is a covering 
$q:\wtl{\Delta}_u\rightarrow\Gamma_w$ with $p=rq$. 
Equivalently, $p$ is universal when any other covering of $\Delta$ is intermediate
to it.

To construct a universal covering, one mimics a standard construction in topology,
%(eg: \cite{Hatcher02}*{\S 1.3}) 
taking as the vertices the homotopy classes of paths in $\Delta$ starting at $v$.
There is an edge
$\wtl{e}$ of $\wtl{\Delta}$ with start vertex the class of $\gamma_1$ and 
finish vertex the class of $\gamma_2$ if and only if there is an edge $e$ of $\Delta$
with $\gamma_1e$ homotopic to $\gamma_2$. 
Define $\wtl{\Delta}_{[v]}\rightarrow\Delta_v$
by sending the class of $\gamma$ to $t(\gamma)$ and the edge $\wtl{e}$ described
above to $e$.

\begin{proposition}
\label{topology:coverings:result600}
$\wtl{\Delta}_{[v]}$ is connected, simply connected, and the map 
$\wtl{\Delta}_{[v]}\rightarrow\Delta_v$ is
a universal covering.
\end{proposition}

\begin{proof}
If $[\gamma]$ is a vertex of $\wtl{\Delta}_{[v]}$ with $\gamma=e_1e_2\ldots e_k$,
then $\wtl{e}_1\wtl{e}_2\ldots\wtl{e}_k$ is a path from $[v]$ to $[\gamma]$
and so $\wtl{\Delta}_{[v]}$ is connected.
That $p$ is a covering is straight forward, and hence $\wtl{\Delta}_{[v]}$
is simply connected, for $\wtl{\gamma}=\wtl{e}_1\wtl{e}_2\ldots\wtl{e}_k$ is 
a closed path at $[v]$ if and only if $e_1e_2\ldots e_k$ is homotopic to
$v$, ie: $\gamma=
e_1e_2\ldots e_k$ is homotopically trivial in $\Delta$, giving $\wtl{\gamma}$
homotopically trivial as $p(\wtl{\gamma})=\gamma$ and by homotopy lifting.
If $r:\Gamma_w\rightarrow\Delta_v$ is a covering then define 
$q:\wtl{\Delta}_{[v]}\rightarrow\Gamma_w$ by $q[\gamma]=t(\gamma')$ where $\gamma'$
is the lift via $r$ of $\gamma$ to the vertex $w$; if $\wtl{e}$ is an edge
with start vertex $[\gamma]$ then let $q(\wtl{e})$ be the lift via $r$ of $p(\wtl{e})$
to the vertex $t(\gamma')$. It is easy to see that $p=rq$ and hence $q$ a covering by
Lemma \ref{topology:coverings:result300}(i).
\qed
\end{proof}

Many authors, anticipating the Galois correspondence, 
define a covering to be universal iff it is simply connected.

\subsection{The lattice of intermediate coverings}
\label{topology:lattice}

Throughout this section $\Lambda,\Delta$ are connected graphs and
$p:\Lambda_u\rightarrow\Delta_v$ is some fixed pointed covering.
A connected pointed intermediate covering 
$\Lambda_u{\rightarrow}\Gamma_x{\rightarrow}\Delta_v$ is
{\em equivalent\/} to 
another such,
$\Lambda_u{\rightarrow}\Upsilon_{y}{\rightarrow}\Delta_v$,
if and only if there is a homeomorphism $\bb:\Gamma_x\rightarrow\Upsilon_{y}$
such that 
$$
\begin{pspicture}(0,0)(6,1.5)
%\showgrid
\rput(3,1.5){$\Gamma_x$}\rput(3,0){$\Upsilon_{y}$}
\psline[linewidth=.1mm]{->}(3,1.2)(3,.4)
\rput(1.8,.75){$\Lambda_u$}\rput(4.2,.75){$\Delta_v$}
\psline[linewidth=.1mm]{->}(2.1,.5)(2.7,.14)\psline[linewidth=.1mm]{->}(2.1,.9)(2.7,1.26)
\psline[linewidth=.1mm]{->}(3.3,.14)(3.9,.5)\psline[linewidth=.1mm]{->}(3.3,1.26)(3.9,.9)
\end{pspicture}
$$
commutes. Let 
$\mathcal{L}(\Lambda_u,\Delta_v)$ be the set of equivalence classes of 
such connected intermediate
coverings. %and we will denote a class by some chosen representative. 

Define $(\Lambda_u{\rightarrow}\Gamma_{x_1}{\rightarrow}\Delta_v)\leq
(\Lambda_u{\rightarrow}\Upsilon_{x_2}{\rightarrow}\Delta_v)$, or just
$\Gamma_{x_1}\leq\Upsilon_{x_2}$, if and only if there is a covering
$s:\Upsilon_{x_2}{\rightarrow}\Gamma_{x_1}$ with $p=r_1sq_2$, where
$r_1$ is the covering $\Gamma_{x_1}\rightarrow\Delta_v$ and $q_2$ is
$\Lambda_u\rightarrow\Upsilon_{x_2}$.
If 
$\bb_1:\Gamma_{x_1}\rightarrow\Psi_{y_1}$ and 
$\bb_2:\Upsilon_{x_2}\rightarrow\Phi_{y_2}$ are homeomorphisms realizing equivalent
coverings, then $\bb_2^{-1}s\bb_1:\Phi_{y_2}\rightarrow\Psi_{y_1}$ is a 
covering with $p=r'_1(\bb_2^{-1}s\bb_1)q'_2$. Thus, $\leq$ is well defined upto 
equivalence, giving $\mathcal{L}(\Lambda_u,\Delta_v)$ the structure of a poset.
We will also write
$\Gamma_w\in\mathcal{L}(\Lambda_u,\Delta_v)$ for an equivalence class of 
intermediate coverings, without reference to the intermediate 
covering maps.

Recall that a poset $(\mathcal{L},\leq)$ equipped with a join $\vee$ (or supremum) and 
meet $\wedge$ (or infimum) is a lattice.
A $\hat{0}$ (resp. $\hat{1}$) is an element such that $\hat{0}\leq x$
(resp. $x\leq\hat{1}$) for all $x$, and 
a lattice isomorphism (resp. anti-isomorphism) 
$\Theta:\mathcal{L}_1\rightarrow\mathcal{L}_2$ is an order-preserving
(resp. order-reversing) bijection whose inverse is also order order-preserving 
(resp. order-reversing). In particular,
an isomorphism sends joins to joins (and meets to meets) and an 
anti-isomorphism sends joins to meets (and meets to
joins). 
A canonical example is
the subgroups of a group $G$, ordered by
inclusion, and with 
$A\vee B=\langle A,B\rangle$, $A\wedge B=A\cap B$, $\hat{0}$ the 
trivial subgroup and $\hat{1}=G$. The remainder of this section is devoted to showing
that $\mathcal{L}(\Lambda_u,\Delta_v)$ is a lattice.

Let $\Lambda_u{\rightarrow}\Gamma_{x_1}{\rightarrow}\Delta_v$ and 
$\Lambda_u{\rightarrow}\Upsilon_{x_2}{\rightarrow}\Delta_v$
be intermediate to $p$, and $q$ the quotient map used in the construction of the pushout 
of the coverings $q_1:\Lambda_u{\rightarrow}\Gamma_{x_1}$
and $q_2:\Lambda_u{\rightarrow}\Upsilon_{x_2}$.
Let $x=q(x_1)=q(x_2)$ and $(\Gamma\coprod_\Lambda\Upsilon)_x$ the resulting
pointed pushout.

\begin{proposition}\label{topological:lattice:result100}
(i). We have the intermediate covering 
$\Lambda_u\stackrel{t_iq_i}{\longrightarrow}
(\Gamma\coprod_\Lambda\Upsilon)_x\stackrel{r}{\rightarrow}\Delta_v$, where $r$ is 
provided by the universality of the pushout. 

\noindent(ii).
Let
$\Psi_{y_1},\Phi_{y_2}\in\mathcal{L}(\Lambda_u,\Delta_v)$ 
be equivalent to $\Gamma_{x_1},\Upsilon_{x_2}$ with
$\bb_1:\Gamma_{x_1}\rightarrow\Psi_{y_1}$,
$\bb_2:\Upsilon_{x_2}\rightarrow\Phi_{y_2}$ the corresponding homeomorphisms and 
$\bb_1\amalg\bb_2:\Gamma\amalg\Upsilon\rightarrow\Psi\amalg\Phi$ 
(disjoint unions) defined by $\bb_1\amalg\bb_2|_{\Gamma}=\bb_1$ and
$\bb_1\amalg\bb_2|_\Upsilon=\bb_2$. Then the map 
$\bb:(\Gamma\coprod_\Lambda\Upsilon)_x\rightarrow(\Psi\coprod_\Lambda\Phi)_{y}$
defined by $\bb q=q'(\bb_1\amalg\bb_2)$ is a homeomorphism making these pointed
pushouts equivalent.
\end{proposition}

Thus there is a well defined pushout of two elements of 
$\mathcal{L}(\Lambda_u,\Delta_v)$.
As the proof will show, the maps 
$t_1,t_2:\Gamma_{x_1},\Upsilon_{x_2}\rightarrow
(\Gamma\coprod_\Lambda\Upsilon)_x$ are coverings, and
so the pushout is a lower bound for 
$\Gamma_{x_1},\Upsilon_{x_2}\in(\mathcal{L}(\Lambda_u,\Delta_v),\leq)$, 
and the universality
implies that it is an infimum.

\begin{proof}[of Proposition \ref{topological:lattice:result100}]
(i). If $v,v'\in\Gamma_1\coprod\Gamma_2$ are vertices with $v\sim v'$ and 
$e$ an edge with start $v$ then, by successively lifting and covering, one
can show that
there is an edge $e'$ with start $v'$ such that $e\sim e'$. 
Now, $v_1$ maps via
$t_1$ to $[v_1]$ and if $[e']$ an edge starting at this vertex
then
$s(e')\sim v_1$, and so by the above there is an edge $e$ starting at 
$v_1$ with $[e']=[e]$. Thus $t_1(e)=[e']$, and so $t_1$ maps the
edges starting at $v_1$ surjectively onto those starting at $t_1(v_1)$.
If $e,e'$ are edges starting at $v_1$ with $t_1(e)=t_1(e')$, then
one gets by induction that $q_1(e)=q_1(e')$, and $q_1$ a covering forces
$e=e'$, and thus $t_1$ (similarly $t_2$) is a covering, hence the $t_iq_i:
\Lambda_u\rightarrow(\Gamma\coprod_\Lambda\Upsilon)_x$ are too. The map $r$
is provided by the universality and is a covering by Lemma 
\ref{topology:coverings:result300}(i). Part (ii) is a tedious but routine diagram
chase.
\qed
\end{proof}

Now to pullbacks. With $\Lambda_u{\rightarrow}\Gamma_{x_1}{\rightarrow}\Delta_v$ and 
$\Lambda_u{\rightarrow}\Upsilon_{x_2}{\rightarrow}\Delta_v$ intermediate to $p$,
$x=x_1\times x_2$ is a vertex of the pullback of the
coverings $r_1:\Gamma_{x_1}{\rightarrow}\Delta_v$ and 
$r_2:\Upsilon_{x_2}{\rightarrow}\Delta_v$. Let 
$(\Gamma\prod_{\Delta}\Upsilon)_x$ be the pointed pullback 
consisting of the component containing the vertex $x$.

\begin{proposition}\label{topological:lattice:result200}
(i). We have the intermediate covering $\Lambda_u\stackrel{q}{\rightarrow}
(\Gamma\prod_\Delta\Upsilon)_x\stackrel{r_it_i}{\longrightarrow}\Delta_v$, 
where $q$ is
provided by the universality of the pullback. 

\noindent(ii). Let
$\Psi_{y_1},\Phi_{y_2}\in\mathcal{L}(\Lambda_u,\Delta_v)$ 
be equivalent to $\Gamma_{x_1},\Upsilon_{x_2}$ with
$\bb_1,\bb_2$ the corresponding homeomorphisms
and 
$\bb:(\Gamma\prod_\Delta\Upsilon)_x\rightarrow(\Psi\prod_\Delta\Phi)_{y}$.
defined by $\bb(x\times y)=\bb_1(x)\times\bb_2(y)$.
Then $\bb$ is a homeomorphism making the pointed
pullbacks equivalent.
\end{proposition}

Thus there is a well defined pullback of two elements of 
$\mathcal{L}(\Lambda_u,\Delta_v)$.
Again the proof shows that the maps 
$t_1,t_2:(\Gamma\prod_\Delta\Upsilon)_x
\rightarrow\Gamma_{x_1},\Upsilon_{x_2}$ 
are coverings, and
so the pullback is an upper bound for 
$\Gamma_{x_1},\Upsilon_{x_2}\in(\mathcal{L}(\Lambda_u,\Delta_v),\leq)$, 
and the universality
implies that it is a supremum.

\begin{proof}[of Proposition \ref{topological:lattice:result200}]
(i). We show that $t_1$ is a covering; $t_2$ is similar. 
From $t_1(e_1\times e_2)=e_1$ it is clear that 
$t_1$ is dimension preserving. For $t_1(u_1\times u_2)=u_1$, 
let $e_1$ be an edge of $\Gamma$
with $s(e_1)=u_1$. Then $e_1$ covers $r_1(e_1)$ which lifts via the 
covering $r_2$ to $u_2$ in $\Upsilon$ to
an edge $e_2$ covering $r_1(e_1)$, ie: with $r_2(e_2)=r_1(e_1)$. 
Thus there is an edge $e_1\times e_2$ of the
pullback with $s(e_1\times e_2)=u_1\times u_2$ and 
$t_1(e_1\times e_2)=e_1$, giving the surjectivity of
$t_1$ on the edges starting at $u_1$. If $e'_1\times e'_2$ starts at 
$u_1\times u_2$ and 
$t_1(e'_1\times e'_2)=e_1$ then $e'_1=e_1$. We have 
$t_2(e'_1\times e'_2)=e'_2$ starting at $u_2$,
and $r_2(e'_2)=r_1t_1(e'_1\times e'_2)=r_1(e_1)=r_2(e_2)$. Thus, as 
$r_2$ is a cover, we have $e'_2=e_2$
and so $e'_1\times e'_2=e_1\times e_2$, and $t_1$ is indeed a covering.
Hence the $r_it_i$ are too and $q$ by Lemma 
\ref{topology:coverings:result300}(i). Part (ii) is an analogous to 
that of Proposition \ref{topological:lattice:result100}.
\qed
\end{proof}

The proof of Proposition \ref{topological:lattice:result200} also shows that
the $t_1,t_2:\Gamma\prod_\Delta\Upsilon\rightarrow\Gamma,\Upsilon$ are 
coverings in the unpointed case.
We pause to observe a slight asymmetry to the duality between pushouts
and pullbacks:
given coverings $r_1,r_2:\Gamma,\Upsilon\rightarrow\Delta$, the
$t_1,t_2:\Gamma\prod_\Delta\Upsilon\rightarrow\Gamma,\Upsilon$ are
coverings, whereas coverings $q_1,q_2:\Lambda\rightarrow\Gamma,\Upsilon$
do not necessarily give coverings $t_1,t_2:\Gamma,\Upsilon\rightarrow
\Gamma\coprod_\Lambda\Upsilon$, unless the $q_i$ are intermediate
$\Lambda\rightarrow(\Gamma\text{ or }\Upsilon)\rightarrow\Delta$.
Indeed, taking the $\Gamma=\Upsilon$ to be two copies of the left hand graph,
$$
\begin{pspicture}(0,0)(13,1)
%\showgrid
\rput(.4,0){
\rput(-.3,.5){$\Gamma=\Upsilon=$}
\rput(1.5,.5){\BoxedEPSF{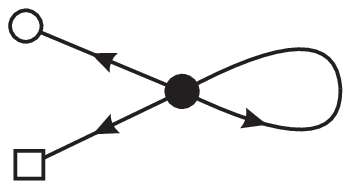 scaled 500}}
\rput(8.45,.5){$=\Lambda=$}
\rput(6,.5){\BoxedEPSF{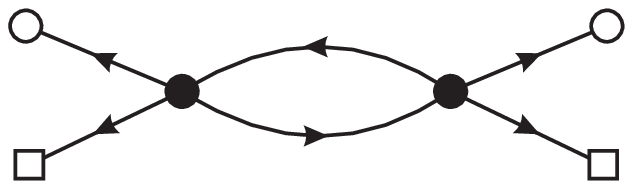 scaled 500}}
\rput(11,.5){\BoxedEPSF{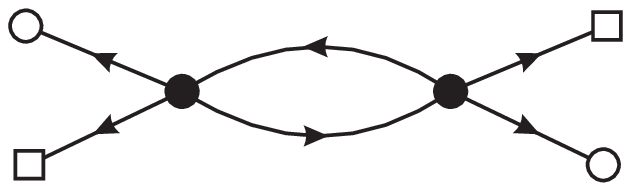 scaled 500}}
\rput(-.7,0){
\psline[linewidth=.1mm]{->}(4.5,.5)(3.5,.5)
\rput(4,.7){$q_i$}
}}
\end{pspicture}
$$
and the coverings $q_i:\Lambda{\rightarrow}\Gamma\text{ or }\Upsilon$ 
(described here by drawing the fibers of the vertices),
then the $t_i$ provided by the pushout construction are not coverings of the pushout.

Summarising the results of this section:

\begin{theorem}[lattice of intermediate coverings]
$\mathcal{L}(\Lambda_u,\Delta_v)$ is a lattice %under $\leq$ 
with join $\Gamma_{x_1}\vee\Upsilon_{x_2}$ the pullback 
$(\Gamma\prod_\Delta\Upsilon)_{x_1\times x_2}$, meet $\Gamma_{x_1}\wedge\Upsilon_{x_2}$
the pushout $(\Gamma\coprod_\Lambda\Upsilon)_{q(x_i)}$, $\widehat{0}=\Delta_v$
and $\widehat{1}=\Lambda_u$.
\end{theorem}

The pointing of the covers in this section is essential if one wishes to work with
{\em connected\/} intermediate coverings and also have a lattice structure 
(both of which
we do). The problem is the pullback: because it is not in general connected, 
we need the
pointing to tell us which component to choose.

\section{The Galois theory of graphs}\label{topology:galoistheory}

The ``Galois correspondence'' between coverings of graphs and subgroups of the
fundamental group goes back to Reidemeister \cite{Reidemeister28}
(see eg: \cite{Collins98}). We provide
a slightly alternative formulation that exploits the lattice structure
of \S \ref{topology:lattice} and is more in the spirit of classical Galois theory.

Throughout this section $p:\Lambda_u\rightarrow\Delta_v$ is a fixed covering
with $\Lambda,\Delta$ connected.
An {\em automorphism\/} (or {\em deck 
transformation\/})
of $p$ is a graph 
homeomorphism $\aa:\Lambda_u\rightarrow\Lambda_{u'}$ making the diagram,
$$
\begin{pspicture}(0,0)(4,1.5)
%\showgrid
\rput(0,-.25){
\rput(0.8,1.5){$\Lambda_u$}
\rput(3.2,1.5){$\Lambda_{u'}$}
\rput(1.95,0.4){$\Delta_v$}
\rput(2,1.7){$\aa$}
\rput(1.195,.825){$p$}
\rput(2.78,.825){$p$}
\psline[linewidth=.1mm]{->}(1.2,1.5)(2.8,1.5)
\psline[linewidth=.1mm]{->}(1,1.25)(1.75,.55)
\psline[linewidth=.1mm]{->}(3,1.25)(2.25,.55)
}
\end{pspicture}
$$
commute. The automorphisms form a group 
$\gal(\Lambda_u,\Delta_v)=\gal(\Lambda_u\stackrel{p}{\rightarrow}\Delta_v)$, 
the {\em Galois group\/} of the covering. 

\begin{lemma}\label{topological:galois:result100}
(i). The action of 
$\gal(\Lambda_u,\Delta_v)$ on $\Lambda$
is orientation preserving.
(ii). The effect of an automorphism $\aa:\Lambda_u\rightarrow\Lambda_{u'}$
is completely determined by $\aa(u)=u'$. In particular, 
the Galois group acts freely
on $\Lambda$.
\end{lemma}

\begin{proof}
(i). Both the edge $e$ and $\aa(e)$ lie in the same fiber of the covering, so that
if $\aa(e)=e^{-1}$ then $p(e)=p(e)^{-1}$, a contradiction, so the Galois group
acts without inversions.
(ii). If $x$ is a vertex of $\Lambda$ and $\gamma$ a path from $u$ to $x$, then
$\aa(x)$ is the terminal vertex of the lift to $u'$
of the path $p(\gamma)$. The images of the edges are handled similarly.
\qed
\end{proof}

The explicit construction of automorphisms is achieved by the following technical
result:

\begin{proposition}\label{topological:galois:result200}
Let $p:\Lambda_{u}\rightarrow\Delta_v$ be a covering
and $u'$ another vertex in the fiber of $v$
such that for any
closed path $\gamma$ at $v$ with lifts $\gamma_i$ 
at $u,u'$, we have $\gamma_1$ closed if and only if 
$\gamma_2$ closed.
%every
%closed path at $v$ lifts to either closed paths at both $u_i$ or
%non-closed paths at both $u_i$. 
%
For any vertex $x\in\Lambda$ and 
path $\mu$ from $u$ to $x$, let $\aa(x)$ be the terminal vertex of
the lift at $u'$ of $p(\mu)$. Then $x\mapsto\aa(x)$ 
extends to an automorphism 
$\aa\in\gal(\Lambda_{u},\Delta_v)$.
\end{proposition}

In particular, for a covering satisfying (i) and (ii) of the Proposition, there is
an element of the Galois group sending the vertex $u$ to the vertex
$u'$. 

\begin{proof}
If $\mu'$ is another path from $u$ to $x$, then 
$p(\mu)p(\mu')^{-1}$ is a closed path at $v$ that lifts to a closed path 
(ie: $\mu(\mu')^{-1}$) at $u$, hence to a 
closed path 
at $u'$. Thus $\aa(x)$ is also the terminal vertex of the lift at $u'$ of 
$p(\mu')$
and $\aa$ is a well defined map $V_\Lambda\rightarrow V_\Lambda$. To extend
$\aa$ to the edges, let $\mu$ be a path from $u$ to the vertex $s(e)$,
and lift the path $p(\mu e)$ to $u'$. Define $\aa(e)$ to be the lift
of $p(e)$ to the terminal vertex of $p(\mu)$. It is easy to see that
$\aa:\Lambda\rightarrow\Lambda$ 
is a surjective dimension-preserving map of graphs, and that $x$ and 
$\aa(x)$ lie in the same fiber of the covering, for any cell
$x$, whence $p\aa=p$. 

It remains to show that $\aa$ is injective. 
For vertices $x$ and $x'$, choose paths $\mu,\mu'$ from $u$ to
$x$ and $x'$. Then if $\aa(x)=\aa(x')$, the lifts at $u'$ of $p(\mu)$
and $p(\mu')$ finish at the same vertex, and so therefore must
$p(\mu)$ and $p(\mu')$, as $p$ is well defined at the vertex
$\aa(x)=\aa(x')$. Thus $\mu,\mu'$ finish
at the same vertex and so $x=x'$. For edges $e,e'$ with $\aa(e)=\aa(e')$,
the injectivity of $\aa$ on the vertices gives that they must have the
same start vertex, and moreover must lie in the same fiber, hence $e=e'$,
by the injectivity of coverings on the edges with start a given vertex.
\qed
\end{proof}

Let
$\Lambda_u{\rightarrow}\Gamma_x{\rightarrow}\Delta_v$ 
be a covering intermediate to $p$ and consider those 
$\aa\in\gal(\Lambda_u,\Delta_v)$
such that 
$$
\begin{pspicture}(0,0)(4,1.5)
%\showgrid
\rput(0,-.25){
\rput(0.8,1.5){$\Lambda_u$}
\rput(3.2,1.5){$\Lambda_{u'}$}
\rput(2,0.4){$\Gamma_x$}
\rput(2,1.7){$\aa$}
\rput(1.195,.825){$q$}
\rput(2.785,.825){$q$}
\psline[linewidth=.1mm]{->}(1.2,1.5)(2.8,1.5)
\psline[linewidth=.1mm]{->}(1,1.25)(1.75,.55)
\psline[linewidth=.1mm]{->}(3,1.25)(2.25,.55)
}
\end{pspicture}
$$
commutes. This gives a subgroup that can be identified with
$\gal(\Lambda_u\stackrel{q}{\rightarrow}\Gamma_x)$. If 
$\Lambda_u{\rightarrow}\Upsilon_{y}{\rightarrow}\Delta_v$ 
is an equivalent covering with homeomorphism $\bb:
\Gamma_x\rightarrow\Upsilon_{y}$,
then 
$$
\begin{pspicture}(0,0)(12,1.5)
%\showgrid
\rput(0,-.25){
\rput(0.8,1.5){$\Lambda_u$}
\rput(3.2,1.5){$\Lambda_{u'}$}
\rput(2,0.4){$\Gamma_x$}
\rput(2,1.7){$\aa$}
\rput(1.195,.825){$q$}
\rput(2.785,.825){$q$}
\psline[linewidth=.1mm]{->}(1.2,1.5)(2.8,1.5)
\psline[linewidth=.1mm]{->}(1,1.25)(1.75,.55)
\psline[linewidth=.1mm]{->}(3,1.25)(2.25,.55)
}
\rput(4.5,.75){commutes}
\rput(5.75,.7){$\Leftrightarrow$}
\rput(5.75,-.25){
\rput(0.8,1.5){$\Lambda_u$}
\rput(3.2,1.5){$\Lambda_{u'}$}
\rput(2,0.4){$\Upsilon_y$}
\rput(2,1.7){$\aa$}
\rput(1.05,.825){$q'$}
\rput(2.785,.825){$q'$}
\psline[linewidth=.1mm]{->}(1.2,1.5)(2.8,1.5)
\psline[linewidth=.1mm]{->}(1,1.25)(1.75,.55)
\psline[linewidth=.1mm]{->}(3,1.25)(2.25,.55)
}
\rput(10,.75){commutes.}
\end{pspicture}
$$
Thus $\gal(\Lambda_u\stackrel{q}{\rightarrow}\Gamma_x)
=\gal(\Lambda_u\stackrel{q'}{\rightarrow}\Upsilon_y)$, and
we can associate in a well defined manner a subgroup of the 
Galois group to an element
of the lattice $\mathcal{L}(\Lambda_u,\Delta_v)$.

On the other hand, if $H\subset\gal(\Lambda_u,\Delta_v)$,
then by Lemma \ref{topological:galois:result100} we may form the 
quotient $\Lambda/H$, and indeed,

\begin{lemma}\label{topological:galois:result400}
If $H_1\subset H_2\subset\gal(\Lambda_u,\Delta_v)$ then,
$$
\Lambda_u\stackrel{q_1}{\rightarrow}(\Lambda/H_1)_{q_1(u)}
\stackrel{s}{\rightarrow}(\Lambda/H_2)_{q_2(u)}
\stackrel{r}{\rightarrow}\Delta_v
$$
are all coverings, where $q_i:\Lambda\rightarrow\Lambda/H_i$ are the
quotient maps and $s,r$ are defined by $sq_1=q_2$ and $rq_2=p$.
\end{lemma}

\begin{proof}
All the complexes are connected, and the $q_i$ coverings by
Lemmas \ref{topological:galois:result100} and
\ref{topology:coverings:result300}(ii).
Thus $s$ is a covering by Lemma \ref{topology:coverings:result300}(i), and
another application, this time to $p=r(sq_1)$, gives that $r$ is a covering.
\qed
\end{proof}

Thus, letting $H_1=H_2=H$, we can associate to
$H\subset\gal(\Lambda_u\,\Delta_v)$ the intermediate covering 
$\Lambda_u{\rightarrow}(\Lambda/H)_{q(u)}{\rightarrow}\Delta_v$,
and by passing to its equivalence class, we get an element of the
lattice $\mathcal{L}(\Lambda_u,\Delta_v)$ associated to $H$.

\begin{proposition}\label{topological:galois:result500}
The following are equivalent for a covering 
$p:\Lambda_u\rightarrow\Delta_v$:
\begin{enumerate}
\item For all closed paths $\gamma$ at
$v$, the lifts of $\gamma$ to each vertex of 
$\text{fib}_{\Lambda\rightarrow\Delta}(v)$
are either all closed or all non-closed;
\item $\gal(\Lambda_u\stackrel{p}{\rightarrow}\Delta_v)$ acts 
regularly on 
$\text{fib}_{\Lambda\rightarrow\Delta}(v)$.
\end{enumerate}
\end{proposition}

In any case, we call the covering $\Lambda_u\rightarrow\Delta_v$ {\em Galois},
with {\em regular\/} a common alternative as the second
part of the Proposition makes clear.
It is clear that if $\Lambda_u\rightarrow\Delta_v$ is Galois then
so is $\Lambda_{u'}\rightarrow\Delta_v$ for any other $u'$ in the
fiber of $v$; if 
$\Lambda_u{\rightarrow}\Gamma_w{\rightarrow}\Delta_v$
is intermediate with 
$\Lambda_u{\rightarrow}\Delta_v$ Galois, then
$\Lambda_u{\rightarrow}\Gamma_w$ is Galois.
%Actually the Galois group of the right hand covering is trivial. Despite
%this, we will see that Galois coverings
%are fairly numerous.

\begin{proof}
The equivalence follows immediately from Proposition 
\ref{topological:galois:result200} and the fact that automorphisms send
closed paths to closed paths and non-closed paths to non-closed paths.
\qed
\end{proof}

If the covering $p:\Lambda_u\rightarrow\Delta_v$ is Galois, 
let $g\in\pi_1(\Delta,v)$ with representative path $\gamma$
and $\aa_g$ an automorphism in
$\gal(\Lambda_u,\Delta_v)$ that sends $u$ to the terminal vertex
$u'$ of the lift of $\gamma$ to $u$. By homotopy lifting, $\aa_g$ depends
only on the vertices $u,u'$ and not on the chosen representative path $\gamma$,
and so the map $\pi_1(\Delta,v)\rightarrow\gal(\Lambda_u,\Delta_v)$
given by $g\mapsto\aa_g$ is well defined. 

\begin{proposition}\label{galois:result500}
If $p:\Lambda_u\rightarrow\Delta_v$ is Galois 
then $g\mapsto \aa_g$ is a surjective homomorphism with kernel 
$p^*\pi_1(\Lambda,u)$, such that under the induced isomorphism
$$
\pi_1(\Delta,v)/p^*\pi_1(\Lambda,u)\rightarrow
\gal(\Lambda_u,\Delta_v),
$$
if
$\Lambda_u\stackrel{q}{\rightarrow}\Gamma_x\stackrel{r}{\rightarrow}\Delta_v$
is intermediate, then the subgroup $r^*\pi_1(\Gamma,x)/p^*\pi_1(\Lambda,u)$
has image $\gal(\Lambda_u\stackrel{q}{\rightarrow}\Gamma_x)$.
\end{proposition}

\begin{proof}
It is easy to check that $\aa_{g_1g_2}=\aa_{g_1}\aa_{g_2}$ and so we have
a homomorphism. If $\bb\in
\gal(\Lambda_u,\Delta_v)$ then $\bb$ is the unique automorphism sending 
$u$ to $\bb(u)$. Taking a path $\mu$ in $\Lambda$ from $u$ to $\bb(u)$
thus gives $\bb=\aa_g$ for $g$ the homotopy class of $p(\mu)$, 
and hence the homomorphism is surjective. 
Because automorphisms act freely, an element $g$ is in the kernel iff
$\aa_g$ fixes the vertex $u$,
and this happens precisely when $g$ can be represented by a path lifting to a 
closed path at $u$, 
ie: when $g\in p^*\pi_1(\Lambda,u)$. It is easy to check that this homomorphism
maps $r^*\pi_1(\Gamma,x)$ onto $\gal(\Lambda_u,\Gamma_x)$.
\qed
\end{proof}

\begin{corollary}\label{galois:result600}
A covering $\Lambda_u\stackrel{p}{\rightarrow}\Delta_v$ is Galois if and only if
$p^*\pi_1(\Lambda,u)$ is a normal subgroup of $\pi_1(\Delta,v)$.
\end{corollary}

\begin{proof}
It remains to show the ``if'' part. 
Let $u'$ be a vertex in the fiber of $v$, $\gamma$ a closed path
at $v$ with lifts $\gamma_1,\gamma_2$ at $u$ and $u'$, and $\mu$ a path
from $u$ to $u'$. Let $g$ and $h$ be the homotopy classes of $\gamma$ and
$p(\mu)$. 
Then $\gamma_1$ is closed iff
$g\in p^*\pi_1(\Lambda,u)$ $\Leftrightarrow$ 
$hgh^{-1}\in p^*\pi_1(\Lambda,u)$ by normality, and
this in turn happens precisely when $\mu\gamma_1\mu^{-1}$ is closed
at $u$, ie: when $\gamma_2$ is closed at $u'$. Thus the covering is
Galois.
\qed
\end{proof}

\begin{proposition}\label{topological:galois:result700}
Let $p:\Lambda_u\rightarrow\Delta_v$ be a Galois covering.
If $H\subset\gal(\Lambda_u,\Delta_v)$
then,
$$
[\gal(\Lambda_u,\Delta_v):H]
=\deg(\Lambda/H_{q(u)}\stackrel{r}{\rightarrow}\Delta_v).
$$
\end{proposition}

\begin{proof}
If a group $G$ acts regularly on a set and $H$ is a subgroup, 
then the number of $H$-orbits
is the index $[G:H]$. 
The result follows as the $H$-orbits on the fiber (via $p$) of $v$
are precisely the vertices of $\Lambda/H$ covering $v$ (via $r$).
\qed
\end{proof}

In particular
the Galois group of a Galois covering has order the degree of the covering.
We have now assembled sufficient machinery to prove,

\begin{theorem}[Galois correspondence]
\label{topological:galois:result800}
Let $\Lambda_u{\rightarrow}\Delta_v$ be a Galois covering
with $\mathcal{L}(\Lambda_u,\Delta_v)$ the lattice of equivalence classes
of intermediate coverings and $\gal(\Lambda_u,\Delta_v)$
the Galois group. Then the map that associates to
$\Lambda_u{\rightarrow}\Gamma_x{\rightarrow}\Delta_v
\in \mathcal{L}(\Lambda_u,\Delta_v)$ the subgroup 
$\gal(\Lambda_u,\Gamma_x)$ is a lattice 
anti-isomorphism from $\mathcal{L}(\Lambda_u,\Delta_v)$ to
the lattice of subgroups of $\gal(\Lambda_u,\Delta_v)$.
Its inverse is the map associating to 
$H\subset\gal(\Lambda_u,\Delta_v)$ the element
$\Lambda_u{\rightarrow}\Lambda/H_{q(u)}
{\rightarrow}\Delta_v\in
\mathcal{L}(\Lambda_u,\Delta_v)$.
\end{theorem}

\begin{proof}
\parshape=9 0pt\hsize 0pt\hsize 0pt\hsize 
0pt.7\hsize 0pt.7\hsize 0pt.7\hsize 0pt.7\hsize 
0pt.7\hsize 0pt\hsize
Let $f$ and $g$ be the two maps described in the theorem. It is easiest to
work from the point of view of $g$: if $H_1\leq H_2$ in the lattice of subgroups,
then the covering $s$ of Lemma \ref{topological:galois:result400}
gives $g(H_2)\leq g(H_1)$, so $g$ is an anti-morphism of lattices.
If $\Lambda_u\rightarrow\Gamma_x\rightarrow\Delta_v$ is intermediate,
then we also have the intermediate covering
$\Lambda_u\rightarrow\Lambda_u/\gal(\Lambda_u,\Gamma_x)\rightarrow\Gamma_x$
with $\Lambda_u\rightarrow\Gamma_x$ is Galois. By Proposition
\ref{topological:galois:result700}, the covering 
$\Lambda_u/\gal(\Lambda_u,\Gamma_x)\rightarrow\Gamma_x$ has degree $1$,
hence is a homeomorphism (Lemma \ref{topological:coverings:result200}(iii))
and we have the diagram at right,
\vadjust{\hfill\smash{\lower -2pt
\llap{
\begin{pspicture}(0,0)(4,1.5)
%\showgrid
\rput(-1,0.2){
\rput(3,1.5){$\Lambda/\gal(\Lambda_u,\Gamma_x)$}\rput(3,0){$\Gamma_{x}$}
\psline[linewidth=.1mm]{->}(3,1.2)(3,.4)
\rput(1.3,.75){$\Lambda_u$}\rput(4.7,.75){$\Delta_v$}
\psline[linewidth=.1mm]{->}(1.6,.5)(2.7,.14)
\rput(-.5,0){\psline[linewidth=.1mm]{->}(2.1,.9)(2.7,1.26)}
\psline[linewidth=.1mm]{->}(3.3,.14)(4.4,.5)
\rput(.5,0){\psline[linewidth=.1mm]{->}(3.3,1.26)(3.9,.9)}
}
\end{pspicture}
}}}\ignorespaces
with the whole square and the left triangle commuting by intermediacy, hence the
right triangle commuting as well. Thus, the intermediate coverings
$\Lambda_u\rightarrow\Gamma_x\rightarrow\Delta_v$ and 
$\Lambda_u\rightarrow\Lambda_u/\gal(\Lambda_u,\Gamma_x)\rightarrow\Delta_v$
are equivalent, and we have $gf=\id$.
If $H\subset\gal(\Lambda_u,\Delta_v)$ and $q:\Lambda\rightarrow\Lambda/H$ the
quotient map,
then $q\aa=q$ for any $\aa\in H$ and
so $H\subset\gal(\Lambda_u,(\Lambda/H)_{q(u)})$ with the covering
$\Lambda_u\rightarrow\Lambda/H_{q(u)}$ intermediate, hence Galois. Proposition 
\ref{topological:galois:result700} gives the index of $H$ in 
$\gal(\Lambda_u,(\Lambda/H)_{q(u)})$ to be the degree of the covering
$\Lambda/H_{q(u)}\rightarrow\Lambda/H_{q(u)}$, ie: $H=
\gal(\Lambda_u,(\Lambda/H)_{q(u)})$,
and we have $fg=\id$.
\qed
\end{proof}

As lattice anti-isomorphisms send joins to meets and meets to joins we
have as an immediate corollary that,

\begin{corollary}\label{topological:galois:result900}
Let $\Lambda_u{\rightarrow}\Delta_v$ be Galois
with $\Lambda_u\rightarrow\Gamma_{x}\rightarrow\Delta_v$ and
$\Lambda_u\rightarrow\Upsilon_{y}\rightarrow\Delta_v$ in the lattice
$\mathcal{L}(\Lambda_u,\Delta_v)$ and $H_1,H_2\subset
\gal(\Lambda_u,\Delta_v)$.
Then,
\begin{align*}
\gal(\Lambda_u,(\Gamma\prod_\Delta\Upsilon)_z)
&=\gal(\Lambda_u,\Gamma_{x})\cap\gal(\Lambda_u,\Upsilon_{y}),\\
\gal(\Lambda_u,(\Gamma\coprod_\Lambda\Upsilon)_z)&=
\langle\gal(\Lambda_u,\Gamma_{x}),\gal(\Lambda_u,\Upsilon_{y})\rangle,
\end{align*}
and the intermediate coverings,
\begin{align*}
\Lambda_u\rightarrow\Lambda/\langle H_1,H_2\rangle_{t(u)}\rightarrow\Delta_v
\text{ and }
\Lambda_u\rightarrow(\Lambda/H_1\coprod_\Lambda \Lambda/H_2)_{qq_i(u)}
\rightarrow\Delta_v,\\
\Lambda_u\rightarrow\Lambda/(H_1\cap H_2)_{t(u)}\rightarrow\Delta_v\text{ and }
\Lambda_u\rightarrow(\Lambda/H_1\prod_\Delta\Lambda/H_2)_w\rightarrow\Delta_v.
\end{align*}
are equivalent (where $t:\Lambda\rightarrow\Lambda/\langle H_1,H_2\rangle$ or
$\Lambda/(H_1\cap H_2)$, $q_i:\Lambda\rightarrow\Lambda/H_i$ and $q$ the quotient
from the pushout).
\end{corollary}

(This result is essentially Theorems 4.3 and 5.5 of \cite{Stallings83}, restated in
our terms.)

The universal cover $\wtl{\Delta}_{u}{\rightarrow}\Delta_v$ ($u=[v]$) is 
Galois by Proposition \ref{topology:coverings:result600} and Corollary
\ref{galois:result600}, and by Proposition \ref{galois:result500}
there is an isomorphism
$$
\varphi:\pi_1(\Delta,v)\rightarrow\gal(\wtl{\Delta}_{u},\Delta_v)
$$
such that if 
$\wtl{\Delta}_{u}{\rightarrow}\Gamma_x
\stackrel{r}{\rightarrow}\Delta_v$ 
is intermediate, then the subgroup
$r^*\pi_1(\Gamma,x)$ has image 
$\gal(\wtl{\Delta}_{u},\Gamma_x)$.
Moreover, two intermediate coverings,
$$
\wtl{\Delta}_{u}{\rightarrow}\Gamma_x
\stackrel{r}{\rightarrow}\Delta_v
\text{ and }
\wtl{\Delta}_{u}{\rightarrow}\Upsilon_{y}\stackrel{r'}{\rightarrow}
\Delta_v,
$$
are equivalent if and only if there is a homeomorphism 
$\bb:\Gamma_x\rightarrow\Upsilon_{y}$ with $r=r'\bb$. We thus obtain the
more familar version of the Galois correspondence, as a special
case of Theorem \ref{topological:galois:result800}:

\begin{corollary}[Galois correspondence for the universal cover]
\label{topological:galois:result1000}
The map that associates to a covering $r:\Gamma_x\rightarrow\Delta_v$
the subgroup $r^*\pi_1(\Gamma,x)$ is a lattice anti-isomorphism from
$\mathcal{L}(\wtl{\Delta}_{u},\Delta_v)$ to the lattice of subgroups
of $\pi_1(\Delta,v)$ that sends Galois covers to normal subgroups.
Its inverse associates to $H\subset\pi_1(\Delta,v)$
the covering $\wtl{\Delta}/\varphi(H)_{q(u)}\rightarrow\Delta_v$.
\end{corollary}

We end this section by showing that the excision of trees has little effect
on the lattice $\mathcal{L}(\Lambda,\Delta)$.
Let $p:\Lambda_u\rightarrow\Delta_v$ be a covering, $T\hookrightarrow\Delta$ a 
spanning tree, $T_i\hookrightarrow\Lambda$ the components of $p^{-1}(T)$ and
$p:(\Lambda/T_i)_{q(u)}\rightarrow(\Delta/T)_{q(v)}$ the induced covering (where
we have (ab)used $q$ for both quotients and $p$ for both coverings).

\begin{theorem}[lattice excision]
\label{galois:lattice:excisetrees}
There is a degree and rank preserving isomorphism of lattices 
$$
\LL(\Lambda,%\stackrel{p}{\rightarrow}
\Delta)\rightarrow
\LL(\Lambda/T_i,%\stackrel{p'}{\rightarrow}
\Delta/T),
$$
that sends the equivalence class of 
$\Lambda_u\rightarrow\Gamma_{x}\stackrel{r}{\rightarrow}\Delta_v$ to the 
equivalence class of 
$\Lambda/T_i\rightarrow\Gamma/T'_i\rightarrow\Delta/T$ 
(with $T'_i\hookrightarrow\Gamma$ the components of $r^{-1}(T)$)
and Galois coverings to Galois coverings.
\end{theorem}

This result could have been shown directly and messily 
at the end of \S \ref{topology:lattice};
we use the Galois correspondence.

\begin{proof}
The quotient $q:\Delta\rightarrow\Delta/T$ induces an isomorphism 
$q^*:\pi_1(\Delta,v)\rightarrow\pi_1(\Delta/T,q(v))$ where $(qp)^*=(pq)^*$ by 
the commutativity of the diagram in Proposition \ref{topology:coverings:result500}.
Thus $q^*p^*\pi_1(\Lambda,v)=p^*\pi_1(\Lambda/T_i,q(u))$ giving that 
$\Lambda\rightarrow\Delta$ is Galois iff $\Lambda/T_i\rightarrow\Delta/T$ is 
Galois by Corollary \ref{galois:result600}, and an isomorphism 
$\gal(\Lambda,\Delta)\rightarrow\gal(\Lambda/T_i,\Delta/T)$ by Proposition
\ref{galois:result500} (leaving off the pointings for clarity). This in turn induces
an isomorphism $\mathcal{L}_1\rightarrow\mathcal{L}_2$ between the subgroup lattices
of these two groups, so that two applications of the Galois correspondence
gives 
$$
\LL(\Lambda,\Delta)\rightarrow\LL_1\rightarrow\LL_2\rightarrow\LL(\Lambda/T_i,\Delta/T),
$$
a composition of an isomorphism and two anti-isomorphisms, hence the result in the
case that $\Lambda\rightarrow\Delta$ is Galois.

If $\Lambda_u\rightarrow\Gamma_{x}\stackrel{r}{\rightarrow}\Delta_v$ is intermediate
then $q^*$ sends the subgroup $r^*\pi_1(\Gamma,x)$ to $r^*\pi_1(\Gamma/T'_i,q(x))$
and so the isomorphism of Galois groups sends $\gal(\Lambda,\Gamma)$ to 
$\gal(\Lambda/T_i,\Gamma/T'_i)$. This gives the desired image of the intermediate
covering, but also, if $\Lambda\rightarrow\Delta$ is not Galois, then 
$\LL(\Lambda,\Delta)$ embeds as a sublattice of $\LL(\wtl{\Delta},\Delta)$, sent to
the sublattice $\LL(\Lambda/T_i,\Delta/T)$ via the result applied to the Galois
covering $\wtl{\Delta}\rightarrow\Delta$.
\qed
\end{proof}

Thus in particular,
there are homeomorphisms
\begin{align*}
(\Gamma\prod_\Delta \Upsilon)/T_k&\rightarrow 
(\Gamma/T_{1j})\prod_{\Delta/T}(\Upsilon/T_{2j}),\\
(\Gamma\coprod_\Lambda \Upsilon)/T_k&\rightarrow 
(\Gamma/T_{1j})\coprod_{\Lambda/T_i}(\Upsilon/T_{2j}).
\end{align*}
where $\Gamma,\Upsilon$ are intermediate to $\Lambda\rightarrow\Delta$,
and the trees $T_{ij},T_{2j},T_k$ are the components of the 
preimages of $T$ via the various coverings.

\section{Graphs of finite rank}
\label{section:invariants}

This section is devoted to a more detailed study of the form of those covering
graphs $\Lambda\rightarrow\Delta$ where $\rk\Lambda<\infty$.

\begin{lemma}\label{ranks100}
Let $\Gamma$ be a connected graph.
%\begin{enumerate}
(i). If $\rk\Gamma<\infty$ and $\Theta$ the trivial graph, 
then $\rk\bigvee_\Theta\Gamma=\rk\Gamma+1$.
(ii). If $\Gamma_1,\Gamma_2$ are connected of finite rank, and 
$\Theta$ finite, then
$$
\rk\biggl(\Gamma_1\bigvee_{\Theta}\Gamma_2\biggr)=|V_\Theta|-1+\sum\rk\Gamma_i.
$$
%\end{enumerate}
\end{lemma}

%(We will require part 4 only in the case $S=\{u\}$, in which case it is 
%obvious anyway).

\begin{proof}
Part (i) follows from the comments at the end
of \S \ref{topological:pushouts}, and (ii)
by induction on $|\Theta|$ and (i).
\qed
\end{proof}

\begin{proposition}\label{finite_rank_characterisation100}
A connected graph $\Lambda$ has finite rank if and only if $\Lambda$  
decomposes as a wedge sum $\Lambda=\Gamma\bigvee_\Theta\Phi$
with $\Gamma$ finite, locally finite, connected, 
$\Theta$ finite, $\Phi$ a forest, and no two vertices of the image of 
$\Theta\hookrightarrow\Phi$ lying in the same component.
%and no two vertices of $\ss(\Delta)$
%lying in the same component.
\end{proposition}

\begin{proof}
If $\Lambda$ has such a decomposition then $\Phi$ necessarily has finitely
many components and finite rankness follows 
by inductively applying Lemma \ref{ranks100}. 
For the converse, %$\Lambda$ is connected by definition. 
fix a basepoint vertex $v$ and spanning 
tree $T$ so that 
there is a finite set $P$ of arcs of $\Lambda$ not
in $T$. Take paths in $T$ from $v$ to the start and terminal vertices of the
arcs of $P$. Let $\Gamma$ be the union of the edges in $P$ and these paths 
and let $\Phi$ the result of removing from $T$ the 
edges in 
the paths.
\qed
\end{proof}

\begin{lemma}\label{invariants:result400}
%\begin{enumerate}
$\Lambda$ connected is of finite rank if and only if for any vertex $v$, the
spine $\widehat{\Lambda}_v$ is finite, locally finite.
%\item The inclusion $\widehat{\Lambda}_v\hookrightarrow\Lambda$ induces
%an isomorphism $\pi_1(\widehat{\Lambda}_v,v)\rightarrow\pi_1(\Lambda,v)$.
%\end{enumerate}
\end{lemma}

\begin{proof}
If $\Lambda$ has finite rank then we have the wedge sum 
decomposition of Proposition \ref{finite_rank_characterisation100}.
If $e$ is an edge not contained in $\Gamma$ and $\gamma$ a closed path
at $v$ containing $e$, then $e$ is contained in a tree component $T$
of $\Phi$. As this component is wedged onto $\Gamma$ at a single vertex,
that part of $\gamma$ contained in $T$ is a closed path, hence
contains a spur. Thus $e$ is 
contained in no closed reduced path at $v$ and the so spine $\widehat{\Lambda}$
is a subgraph of $\Gamma$, hence finite, locally finite. Conversely, a finite spine has
finite rank, hence so does $\Lambda$ by Lemma \ref{topological:spines:result100}(i).
\qed
\end{proof}

\begin{proposition}%[spine decomposition]
\label{finite_rank_characterisation200}
%A connected graph $\Lambda$ decomposes as a wedge sum,
Let $\Lambda$ be a connected graph, $\Gamma\hookrightarrow\Lambda$ a 
connected subgraph and $v\in\Gamma$ a vertex such that every closed reduced
path at $v$ in $\Lambda$ is contained in $\Gamma$. Then $\Lambda$ has a 
wedge sum decomposition $\Lambda=\Gamma\bigvee_\Theta\Phi$
with $\Phi$ a forest 
and no two vertices of the image of 
$\Theta\hookrightarrow\Phi$ lying in the same component.
\end{proposition}

\begin{proof}
Consider an edge $e$ of $\Lambda\setminus\Gamma$ having 
at least 
one of its end vertices $s(e)$ or
$t(e)$, in $\Gamma$. 
For definiteness we can assume, by relabeling
the edges in the arc $\ov{e}$, that it is $s(e)$ that is a vertex of 
$\Gamma$. If $t(e)\in\Gamma$ 
then by traversing a reduced path in $\Gamma$ from $v$ to $s(e)$, crossing
$e$ and a reduced path in $\Gamma$ from $t(e)$ to $v$, we get a closed 
reduced path not contained in $\Gamma$, a contradiction. 
Thus $t(e)\not\in\Gamma$. Let $T_e$ be the union 
of all the reduced paths
in $\Lambda\setminus\{e\}$ starting at $t(e)$, so we have the situation
as in (a):
$$
\begin{pspicture}(0,0)(13,3)
%\showgrid
\rput(1.5,1.7){\BoxedEPSF{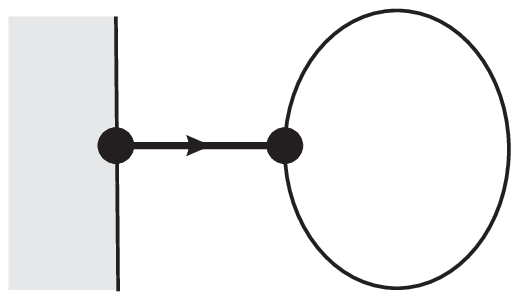 scaled 650}}
\rput(0.1,1.7){$\Gamma$}
\rput(2.35,1.7){$T_e$}\rput(1.1,1.95){$e$}
\rput(1.5,0){(a)}
\rput(6.5,1.7){\BoxedEPSF{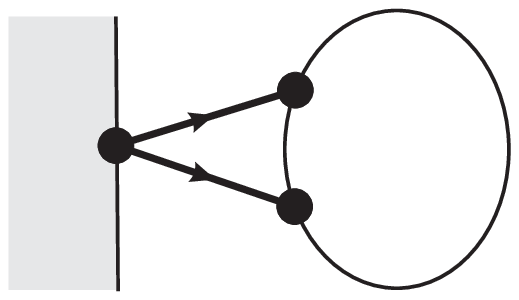 scaled 650}}
\rput(5.1,1.7){$\Gamma$}
\rput(7.4,1.7){$T_e$}
\rput(6.1,2.1){$e$}\rput(6.1,1.3){$e'$}
\rput(6.5,0){(b)}
\rput(11.5,1.7){\BoxedEPSF{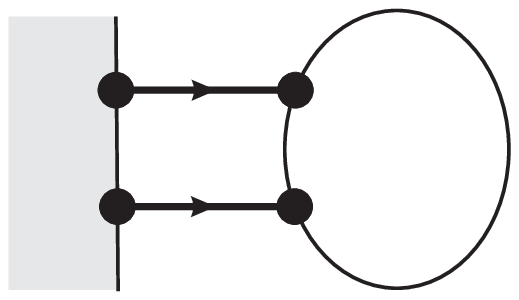 scaled 650}}
\rput(10.1,1.7){$\Gamma$}
\rput(12.4,1.7){$T_e$}
\rput(11.1,2.3){$e$}\rput(11.1,1.1){$e'$}
\rput(11.5,0){(c)}
\end{pspicture}
$$
%Now, $T_e$ is a tree, for 
If $\gamma$ is a non-trivial closed path in
$T_e$ starting at $t(e)$, then a path from $v$ to $t(e)$, 
traversing $\gamma$, and going the same way back to
$v$ cannot be reduced. But the only place a spur can occur is in $\gamma$
and so $T_e$ is a tree.
If $e'$ is another edge of $\Lambda\setminus\Gamma$ with
$s(e')\in\Gamma$ then we claim that neither of the 
two situations (b) and (c) above can occur, ie: $t(e')$ is not a vertex of 
$T_e$. For otherwise, a reduced closed path in $T_e$ from $t(e)$ to
$t(e')$ will give a reduced closed path at $v$ not in $\Gamma$.
Thus, another edge $e'$ yields a tree $T_{e'}$ defined like
$T_e$, but disjoint from it. Each component of $\Phi$ is thus obtained 
this way.
\qed
\end{proof}

In particular we have such a decomposition involving a spine, and
so $\Lambda$ is made up of its spine at some vertex, 
together with a collection of
trees, each connected to $\widehat{\Lambda}_v$ by a single edge. If
$\Lambda$ has finite rank then $\widehat{\Lambda}$ and $\Theta$
are finite, and we have
\begin{equation}\label{finite:rank200}
\Lambda=\biggl(\cdots\biggl(\biggl(\widehat{\Lambda}_v
\bigvee_{\Theta_1} T_1\biggr)
\bigvee_{\Theta_2} T_2\biggr)
\cdots\biggr)\bigvee_{\Theta_k} T_k\biggr),
\end{equation}
with the $\Theta_i$ single vertices, the
$\Theta_i\hookrightarrow\widehat{\Lambda}_v$,
and the images $\Theta_i\hookrightarrow T_i$ having valency one.
Moreover, if $\Lambda\rightarrow\Delta$ is a covering with $\Delta$
single vertexed and $\Lambda$ of finite rank, 
then by Proposition \ref{topology:coverings:result400}(i), each 
tree $T_i$ realizes an embedding $\RR\hookrightarrow\Lambda$ 
of the real line in $\Lambda$, and as the
spine is finite, the trees are thus paired
$$
\begin{pspicture}(0,0)(14,2)
%\showgrid
\rput(3,0){
\rput(4,1){\BoxedEPSF{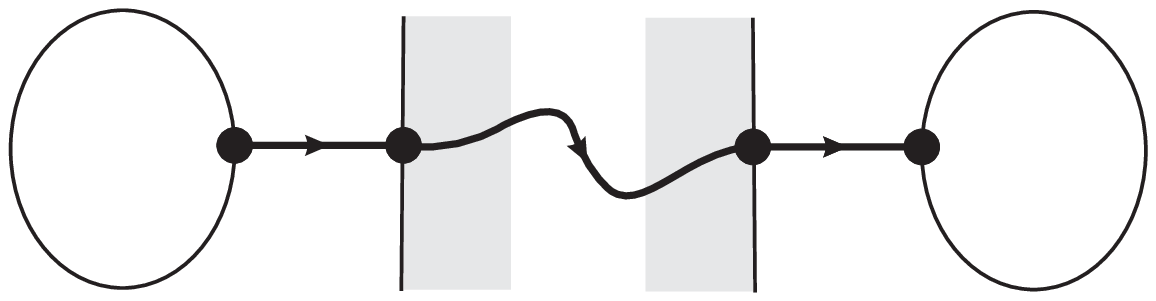 scaled 650}}
\rput(4,1.7){$\widehat{\Lambda}$}
\rput(3.75,.9){$\gamma$}
\rput(2.2,1.3){$e_1$}\rput(5.6,1.3){$e_2$}
\rput(1,1){$T_{e_1}$}\rput(7,1){$T_{e_2}$}
}
\rput(14,1){$(\ddag)$}
\end{pspicture}
$$
with the $e_i$ (and indeed all the edges in the path 
$\mathcal{R}\hookrightarrow\Lambda$) 
in the same fiber of the covering. This pairing will play an important role
in \S \ref{section:pullbacks}.

\begin{corollary}\label{finiterank:result500}
Let $\Lambda\rightarrow\Delta$ be a covering with $\Delta$ non-trivial,
single vertexed and $\rk\Lambda<\infty$. Then 
$\deg(\Lambda\rightarrow\Delta)<\infty$ if and only if $\Lambda=\widehat{\Lambda}_v$.
\end{corollary}

\begin{proof}
If $\Lambda$ is more than $\widehat{\Lambda}_v$ then one of the trees $T_i$ in 
the decomposition (\ref{finite:rank200}) is non trivial and by Proposition
\ref{topology:coverings:result400}(i) we get the real line
$\RR\hookrightarrow\Lambda$, with image in the fiber of an edge, contradicting
the finiteness of the degree. The converse follows from Lemma 
\ref{invariants:result400}.
\qed
\end{proof}

\begin{proposition}
\label{finiterank:result600}
Let $\Lambda\rightarrow\Delta$ be a covering with (i). $\rk\Delta>1$, (ii). 
$\rk\Lambda<\infty$, and (iii). for any intermediate covering 
$\Lambda\rightarrow\Gamma\rightarrow\Delta$ we have $\rk\Gamma<\infty$. Then 
$\deg(\Lambda\rightarrow\Delta)<\infty$.
\end{proposition}

The covering $\RR\rightarrow\Delta$ of a single vertexed $\Delta$
of rank $1$ by the real line shows why the $\rk\Delta>1$ condition 
cannot be dropped.

\begin{proof}
By lattice excision, Theorem \ref{galois:lattice:excisetrees}, we may pass to
the $\Delta$ single vertexed case while preserving (i)-(iii). Establishing
the degree here and passing back to the general $\Delta$ will give the result.
If the degree of the covering
$\Lambda\rightarrow\Delta$ is infinite for $\Delta$ single vertexed, then
by Corollary \ref{finiterank:result500},
in the decomposition (\ref{finite:rank200}) for $\Lambda$, one
of the trees is non-empty and $\Lambda$ has the form
of the graph in Proposition \ref{topology:coverings:result400} with this non-empty
tree the union of the edge $e$ and $\Upsilon_2$. 

Let $\Gamma$ be a graph as defined as follows:
take the union of $\Upsilon_1$, the edge $e$ and
$\aa(\RR)\cap\Upsilon_2$, where $\aa(\mathcal{R})$ is the embedding
of the real line given by Proposition \ref{topology:coverings:result400}(i). 
At each vertex of $\aa(\RR)\cap\Upsilon_2$
place $\rk\Delta-1$ edge loops:
$$
\begin{pspicture}(0,0)(14,2)
%\showgrid
\rput(2,0){
\rput(2,1){\BoxedEPSF{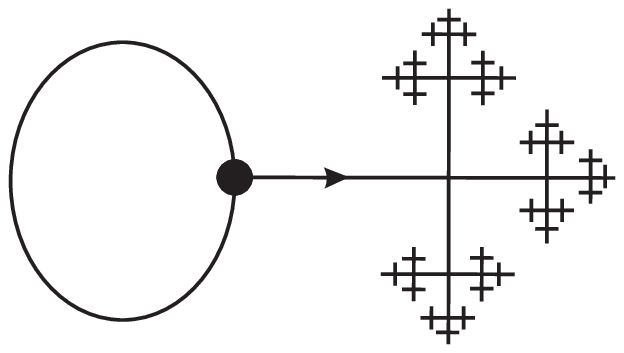 scaled 650}}
\rput(8,1){\BoxedEPSF{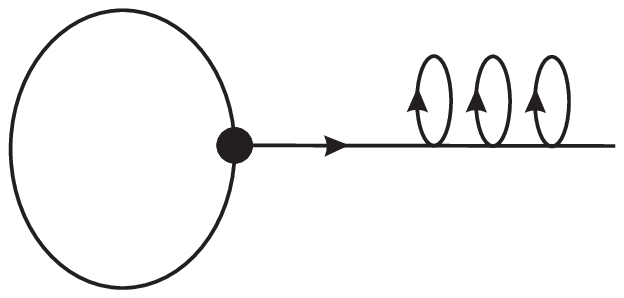 scaled 650}}
\rput(.8,1){$\Upsilon_1$}\rput(6.75,1){$\Upsilon_1$}
\rput(-.2,1.5){$\Lambda$}\rput(5.8,1.5){$\Gamma$}
\rput(9,.7){$\aa(\mathcal{R})\cap\Upsilon_2$}
\rput(10.2,1.3){$\cdots$}
}
\end{pspicture}
$$
(this picture depicting the $\rk\Delta=2$ case).
Then there is an obvious covering $\Gamma\rightarrow\Delta$ so that
by Proposition \ref{topology:coverings:result400}(ii) we have an intermediate covering
$\Lambda\rightarrow\Gamma\rightarrow\Delta$.
Equally obviously, $\Gamma$ has infinite rank, contradicting (iii). Thus, 
$\deg(\Lambda\rightarrow\Delta)<\infty$.
\qed
\end{proof}

\begin{proposition}
\label{finiterank:result700}
%Let $\Psi\rightarrow\Lambda\rightarrow\Delta$ be coverings with
%$\rk\Lambda<\infty$,
%$\deg(\Lambda\rightarrow\Delta)$ infinite, and 
%$\Psi\rightarrow\Delta$ Galois. Then $\Psi$ is simply connected.
Let $\Psi\rightarrow\Lambda\rightarrow\Delta$
be coverings with $\rk\Lambda<\infty$,
$\Psi\rightarrow\Delta$ Galois, and $\Psi$ not simply connected. Then
$\deg(\Lambda\rightarrow\Delta)<\infty$.
\end{proposition}

The idea of the proof is that if the degree is infinite, then
$\Lambda$ has a hanging tree in its 
spine decomposition, and so $\Psi$ does too. But $\Psi$ should look the same
at every point, hence {\em is\/} a tree.

\begin{proof}
Apply lattice excision to 
$\LL(\Psi,\Delta)$, and as $\pi_1(\Psi,u)$ is unaffected by the excision
of trees, we may assume that $\Delta$ is single vertexed.
As $\deg(\Lambda\rightarrow\Delta)$
is infinite, the spine decomposition for $\Lambda$ has an infinite tree,
and $\Lambda$ has the form of Proposition \ref{topology:coverings:result400}. Thus
$\Psi$ does too, by part (iii) of this Proposition, with subgraphs $\Upsilon'_i
\hookrightarrow\Psi$, edge $e'$ and $\Psi'_1$ a tree.
Take a closed
reduced path $\gamma$ in $\Upsilon'_2$, and choose a vertex $u_1$ of $\Upsilon'_1$ 
such that the reduced path from $u_1$ to $s(e')$ has at least as many edges
as $\gamma$. Project $\gamma$ via the covering $\Psi\rightarrow
\Delta$ to a closed reduced path, and then lift to $u_1$. The result is
reduced, closed by Proposition \ref{topological:galois:result500}, and 
entirely contained in the tree $\Upsilon'_1$, hence trivial.
Thus $\gamma$ is also trivial so that
$\Upsilon'_2$ is a tree and $\Psi$ is simply connected. 
\qed
\end{proof}

\begin{proposition}%[topological Marshall Hall]
\label{finiterank:result800}
Let $\Lambda_u\rightarrow\Delta_v$ be a covering with $\rk\Lambda<\infty$
and $\gamma$ a non-trivial reduced closed path at $v$ lifting to a non-closed path 
at $u$. Then there is an intermediate covering 
$\Lambda_u\rightarrow\Gamma_w\rightarrow\Delta_v$ 
with $\deg(\Gamma\rightarrow\Delta)$
finite and $\gamma$ lifting to a non-closed path at $w$.
\end{proposition}

Stallings shows something very similar \cite{Stallings83}*{Theorem 6.1}
starting from a finite immersion rather than a covering. As the proof shows,
the path $\gamma$ in Proposition \ref{finiterank:result800} can be replaced
by finitely many such paths. Moreover, the intermediate $\Gamma$ constructed
has the property that any set of Schreier generators for $\pi_1(\Lambda,u)$
can be extended to a set of Schreier generators for $\pi_1(\Gamma,w)$.

\begin{proof}
If $T\hookrightarrow\Delta$ is a spanning tree and $q:\Delta\rightarrow\Delta/T$
then $\gamma$ cannot be contained in $T$, and so $q(\gamma)$ is non-trivial, closed
and reduced. If the lift of $q(\gamma)$ to $\Lambda/T_i$ is closed then the lift
of $\gamma$ to $\Lambda$ has start and finish vertices that lie in the same
component $T_i$ of $p^{-1}(T)$, mapped homeomorphically onto $T$ by the covering, 
and thus implying that
$\gamma$ is not closed. Thus we may apply lattice excision and pass to the
single vertexed case while maintaining $\gamma$ and its properties. Moreover,
the conclusion in this case gives the result in general as closed paths go to closed paths
when excising trees.

If the lift $\gamma_1$ of $\gamma$
at $u$ is not contained in the spine $\widehat{\Lambda}_u$, then its terminal 
vertex lies in a tree $T_{e_i}$ of the spine decomposition $(\ddag)$. By adding
an edge if necessary to $\widehat{\Lambda}_u\cup\gamma_1$, we obtain
a finite subgraph whose coboundary edges are paired, with the edges in each
pair covering the same edge in $\Delta$, as below left:
$$
\begin{pspicture}(0,0)(14,1)
%\showgrid
\rput(4.05,.5){$\widehat{\Lambda}_u\cup\gamma_1$}
\rput(10,.2){$\Gamma$}
\rput(4,.5){\BoxedEPSF{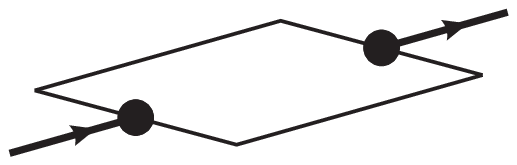 scaled 740}}
\rput(10,.5){\BoxedEPSF{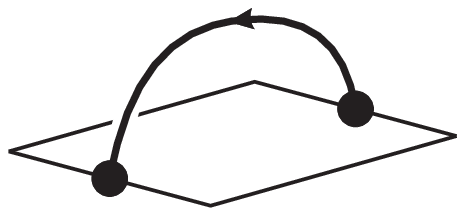 scaled 740}}
\end{pspicture}
$$
(if the lift is contained in the spine,
take $\widehat{\Lambda}_u$ itself).
In any case, let $\Gamma$ be $\widehat{\Lambda}_u\cup\gamma_1$ together with a single 
edge replacing each pair as above right. Restricting the covering 
$\Lambda\rightarrow\Delta$ to $\widehat{\Lambda}_u\cup\gamma_1$ and 
mapping the new edges to the
common image of the old edge pairs gives a finite covering $\Gamma\rightarrow\Delta$
and hence by Lemma \ref{topology:coverings:result300}(i) an intermediate covering
$\Lambda\stackrel{q}{\rightarrow}\Gamma{\rightarrow}\Delta$, with 
$q(\gamma_1)$ non-closed at $q(u)$.
\qed
\end{proof}

For the rest of this section we investigate the rank implications 
of the decomposition (\ref{finite:rank200})
and the pairing $(\ddag)$ in a special case. 
Suppose $\Lambda\rightarrow\Delta$ is a covering 
with $\Delta$ single vertexed, $\rk\Delta=2$, $\Lambda$ non-simply 
connected and $\rk\Lambda<\infty$. Let $x_i^{\pm 1}, (1\leq i\leq 2)$
be the edge loops of $\Delta$
and fix a spine so we have the decomposition (\ref{finite:rank200}).

An {\em extended spine\/} for such a $\Lambda$ 
is a connected subgraph $\Gamma\hookrightarrow\Lambda$ 
obtained by adding finitely many edges to a spine, so that every
vertex of $\Gamma$ is incident with either zero or three edges
in its coboundary $\delta\Gamma$. It is always possible to find
an extended spine: take the 
union of the spine $\widehat{\Lambda}_u$ and each 
edge $e\in\delta\widehat{\Lambda}_u$ in its coboundary. 
Observe that $\Gamma$ is finite and the decomposition
(\ref{finite:rank200}) gives $\rk\Gamma=\rk\widehat{\Lambda}_u=\rk\Lambda$.
Call a vertex of the extended spine $\Gamma$
{\em interior\/} (respectively {\em boundary\/}) when it is incident
with zero (resp. three) edges in $\delta\Gamma$.

We have the pairing of trees $(\ddag)$ for an extended spine, so that
each boundary vertex $v_1$ is paired with another $v_2$,
$$
\begin{pspicture}(0,0)(14,2)
%\showgrid
\rput(3,0){
\rput(4,1){\BoxedEPSF{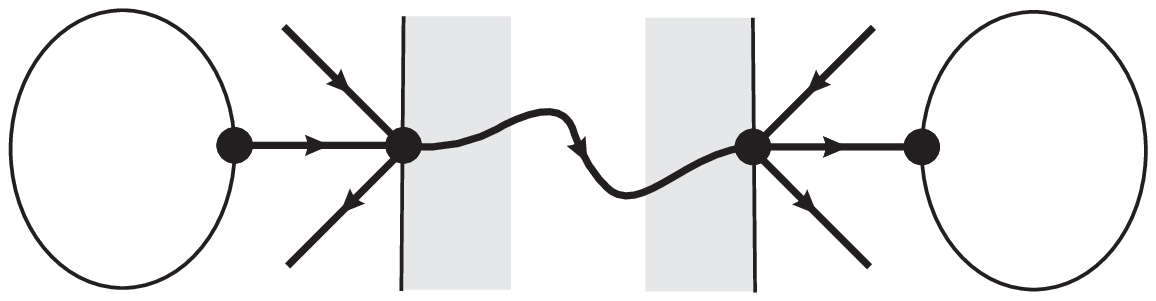 scaled 650}}
\rput(4,1.7){$\Gamma$}
\rput(3.8,.9){$\gamma$}
\rput(2.2,1.25){$e_1$}\rput(5.8,1.25){$e_2$}
\rput(3.1,1.25){$v_1$}\rput(4.9,1.25){$v_2$}
\rput(1,1){$T_{e_1}$}\rput(7,1){$T_{e_2}$}
}
\rput(14,1){$(*)$}
\end{pspicture}
$$
with $e_1,e_2$ and all the edges in the path $\gamma=\aa(\RR)\cap\Gamma$ 
covering an edge loop $x_i\in\Delta$. Call this an {\em $x_i$-pair\/},
($i=1,2$).

For two $x_i$-pairs (fixed $i$), the
respective $\gamma$ paths share no vertices in common, for otherwise
there would be two distinct edges covering the same $x_i\in\Delta$ starting
at such a common vertex. 
Moreover, $\gamma$ must contain vertices of $\Gamma$ apart from
the two boundary vertices $v_1,v_2$, otherwise $\Lambda$ 
would be simply connected.
These other vertices are incident with at least two edges of $\gamma\in\Gamma$, 
hence at most $2$ edges of the coboundary $\delta\Gamma$, and thus must be 
interior.

\begin{lemma}\label{finiterank:result1000}
If $n_i, (i=1,2)$, is the number of $x_i$-pairs in an extended spine $\Gamma$, 
then then number of interior vertices is at least $\sum n_i$.
\end{lemma}

(Lemma \ref{finiterank:result1000} is not true in the 
case $\rk\Delta>2$).

\begin{proof}
The number of interior vertices is $|V_\Gamma|-2\sum n_i$ and the number of
edges of $\Gamma$ is $4(|V_\Gamma|-2\sum n_i)+2\sum n_i$, hence 
$\rk\Gamma-1=|V_\Gamma|-3\sum n_i$ by \S \ref{topological:rank}. 
As $\Lambda$ is not simply connected,
$\rk\Lambda-1=\rk\Gamma-1\geq 0$, thus $|V_\Gamma|-2\sum n_i\geq
\sum n_i$ as required.
\qed
\end{proof}

It will be helpful in \S \ref{section:pullbacks} 
to have a pictorial description of the quantity
$\rk-1$ for our graphs. To this end, a {\em checker\/} is a small plastic
disk, as used in the eponymous boardgame (called {\em draughts\/}
in British English). We place black checkers on some of the
vertices of an extended spine $\Gamma$ according to the 
following scheme: place black checkers on all the interior vertices
of $\Gamma$; for each $x_1$-pair in (*), take the interior vertex
on the path $\gamma$ that is closest to $v_1$ (ie: is the terminal vertex of the
edge of $\gamma$ whose start vertex is $v_1$) and {\em remove\/} its
checker; for each
$x_2$-pair, we can find, by Lemma \ref{finiterank:result1000}, an interior
vertex with a checker still on it. Choose such a vertex and remove its checker
also.

\begin{lemma}\label{finiterank:whitevertices}
With black checkers placed on the vertices of an extended spine for $\Lambda$ 
as above, the number of black checkers is $\rk\Lambda-1$.
\end{lemma}

\begin{proof}
We saw in the proof of Lemma \ref{finiterank:result1000} that $\rk\Lambda-1=
\rk\Gamma-1$ is equal to the number of interior vertices of $\Gamma$ less
the number of $x_i$-pairs $(i=1,2)$.
\qed
\end{proof}

%We remark that the reason we have placed checkers rather than, say, colouring the
%vertices, is that checkers can be easily {\em moved\/}.

From now on we will only use the extended spine obtained by adding the 
coboundary edges to some fixed spine $\widehat{\Lambda}_u$.

Let $p:\Lambda_u\rightarrow\Delta_v$ be a covering with $\rk\Delta=2$,
$\rk\Lambda<\infty$ and $\Lambda$ not simply connected. A spanning tree
$T\hookrightarrow\Delta$ induces a covering $\Lambda/T_i\rightarrow\Delta/T$
with $\Delta/T$ single vertexed. Let $\HH(\Lambda_u\rightarrow\Delta_v)$ be the
number of vertices of the spine of $\Lambda/T_i$ at $q(u)$ and 
$n_i(\Lambda_u\rightarrow\Delta_v)$ the number of $x_i$-pairs in the 
extended spine. The homeomorphism class of $\Lambda/T_i$ and the spine
are independent of the spaning tree $T$, hence the quantities
$\HH(\Lambda_u\rightarrow\Delta_v)$ and $n_i(\Lambda_u\rightarrow\Delta_v)$
are too.

\section{Pullbacks}\label{section:pullbacks}

Let $p_i:\Lambda_i:=\Lambda_{u_i}\rightarrow\Delta_v, (i=1,2)$ be coverings and 
$(\Lambda_{1}\prod_\Delta\Lambda_{2})$ their (unpointed) pullback. 
If $\widehat{\Lambda}_{u_i}$ is the spine at $u_i$ then we can restrict
the coverings to maps $p_i:\widehat{\Lambda}_{u_i}\rightarrow\Delta_v$
and form the pullback $\widehat{\Lambda}_{u_1}\prod_\Delta\widehat{\Lambda}_{u_2}$.

\begin{proposition}[spine decomposition of pullbacks]
\label{pullbacks:spinedecomposition}
The pullback $\Lambda=(\Lambda_{1}\prod_\Delta\Lambda_{2})$ 
has a wedge sum
decomposition 
$\Lambda=(\widehat{\Lambda}_{u_1}\prod_\Delta\widehat{\Lambda}_{u_2})
\bigvee_\Theta\Phi$
with $\Phi$ a forest 
and no two vertices of the image of 
$\Theta\hookrightarrow\Phi$ lying in the same component.
\end{proposition}

\begin{proof}
Let $\Lambda_i=\widehat{\Lambda}_{u_i}\bigvee_{\Theta_i}\Phi_i, (i=1,2)$ be the 
spine 
decomposition, $t_i:\Lambda_1\prod_\Delta\Lambda_2\rightarrow\Lambda_i, (i=1,2)$
the coverings provided by the pullback and
$\Omega$ a connected component of the pullback. If 
$\Omega\cap(\widehat{\Lambda}_{u_1}\prod_\Delta\widehat{\Lambda}_{u_2})
=\varnothing$, then a reduced closed path $\gamma\in\Omega$ must map via one of 
the $t_i$ to a closed path in the forest $\Phi_i$. As the images under 
coverings of
reduced paths are reduced, $t_i(\gamma)$ must contain a spur which can be lifted
to a spur in $\gamma$. Thus $\Omega$ is a tree. 

Otherwise choose 
a vertex $w_1\times w_2$ in 
$\Omega\cap(\widehat{\Lambda}_{u_1}\prod_\Delta\widehat{\Lambda}_{u_2})$ and
let $\Gamma$ be the connected component of this intersection containing
$w_1\times w_2$. If $\gamma$ a reduced closed path at $w_1\times w_2$ then
$t_i(\gamma), (i=1,2)$ a reduced closed path at 
$w_i\in\widehat{\Lambda}_{u_i}$, hence 
by Lemma \ref{topological:spines:result100}(ii), 
$t_i(\gamma)\in\widehat{\Lambda}_{u_i}$ and thus 
$\gamma\in\widehat{\Lambda}_{u_1}\prod_\Delta\widehat{\Lambda}_{u_2}$.
Applying Proposition \ref{finite_rank_characterisation200}, we have $\Omega$
a wedge sum of $\Gamma$ and a forest of the required form.
\qed
\end{proof}

%Indeed, the proof shows that the spine of the pullback is a subgraph of the 
%pullback of the spines. In particular we have,

\begin{corollary}[Howsen-Stallings]
\label{pullbacks:result200}
Let $p_i:\Lambda_i\rightarrow\Delta, (i=1,2),$ be coverings with 
$\rk\Lambda_i<\infty$ and $u_1\times u_2$ a vertex of their pullback. Then 
$\rk(\Lambda_1\prod_\Delta\Lambda_2)_{u_1\times u_2}<\infty.$
\end{corollary}

\begin{proof}
The component $\Omega$ of the pullback containing $u_1\times u_2$ is either a tree
or the wedge sum of a finite graph and a forest as described in Proposition
\ref{pullbacks:spinedecomposition}. Either case gives the result.
\qed
\end{proof}

The remainder of this section is devoted to a proof of an estimate for the
rank of the pullback of finite rank graphs in a special case. Let 
$p_j:\Lambda_j:=\Lambda_{u_j}\rightarrow\Delta_v, (j=1,2)$ be coverings
with $\rk\Delta=2$, $\rk\Lambda_j<\infty$ and the
$\Lambda_j$ not simply connected. 
Let $\HH_j:=\HH(\Lambda_{u_j}\rightarrow\Delta_v)$ and 
$n_{ji}:=n_i(\Lambda_{u_j}\rightarrow\Delta_v)$ be as at the end of 
\S \ref{section:invariants}.

\begin{theorem}\label{pullback:rankestimate}
For $i=1,2$,
$$
\sum_\Omega (\rk\Omega-1)\leq \prod_j(\rk\Lambda_j-1)
%+\varepsilon_{2j}\HH_1+\varepsilon_{1j}\HH_2-\varepsilon_{1j}\varepsilon_{2j},
+\HH_1\HH_2-(\HH_1-n_{1i})(\HH_2-n_{2i}),
$$
the sum over all non simply connected components $\Omega$ of the pullback
$\Lambda_1\prod_\Delta\Lambda_2$.
\end{theorem}

\begin{proof}%[of Theorem \ref{pullback:rankestimate}]
Lattice excision and the definition of the $\HH_j$ and $n_{ji}$ allow us to pass
to the $\Delta$ single vertexed case. 
Suppose then that $\Delta$ has edge loops 
$x_i^{\pm 1}, (1\leq i\leq 2)$ at the vertex $v$,
extended spines 
$\widehat{\Lambda}_{u_j}\hookrightarrow\Gamma_j\hookrightarrow\Lambda_j$, 
and by restricting the covering maps $p_j$ appropriately, the
pullbacks 
$\widehat{\Lambda}_{u_1}\prod_\Delta\widehat{\Lambda}_{u_2}\hookrightarrow
\Gamma_1\prod_\Delta\Gamma_2\hookrightarrow
\Lambda_1\prod_\Delta\Lambda_2$ with $t_j:\Lambda_1\prod_\Delta\Lambda_2
\rightarrow\Lambda_j$ the resulting covering maps. 

Place black checkers on the vertices of the extended spines $\Gamma_j$ as in 
\S \ref{section:invariants} and place a black checker on
a vertex $v_1\times v_2$ of $\Gamma_1\prod_\Delta\Gamma_2$ 
precisely when both $t_j(v_j)\in\Gamma_j, (j=1,2)$ have black checkers
on them.
By Lemma \ref{finiterank:whitevertices}, and the construction of the
pullback for $\Delta$ single vertexed, we get the number of vertices
in $\Gamma_1\prod_\Delta\Gamma_2$ with black checkers is equal to 
$\prod(\rk\Lambda_j-1)$.

Let $\Omega$ be a non simply connected component of the pullback 
$\Lambda_1\prod_\Delta\Lambda_2$ and 
$\Upsilon=\Omega\cap(\Gamma_1\prod_\Delta\Gamma_2)$.
If $v_1\times v_2$ is the start vertex of
at least one edge in the coboundary $\delta\Upsilon$,
then at least one of the $v_j$ must be incident with at least one,
hence three, edges of the coboundary $\delta\Gamma_j$. 
Lifting these three via
the covering $t_j$ to $v_1\times v_2$ gives at least three edges 
starting at $v_1\times v_2$ in the coboundary
$\delta\Upsilon$. Four coboundary edges starting here
would mean that $\Omega$ was simply connected,
hence every vertex of $\Upsilon$
is incident with either
zero or three coboundary edges. 

We can thus extend the interior/boundary terminology of \S \ref{section:invariants}
to the vertices of 
$\Upsilon$, and observe that a vertex of $\Upsilon$ covering,
via either of the $t_j$, a boundary vertex $\in\Gamma_j$, must itself be a
boundary vertex.
The upshot is that $\Upsilon$ is an extended spine 
in $\Omega$ and by Proposition \ref{pullbacks:spinedecomposition}, $\rk\Omega-1=
\rk\Upsilon-1$. Now place {\em red\/} checkers on the vertices of $\Upsilon$  
as in \S \ref{section:invariants} and do this for each non-simply connected 
component $\Omega$. The number of red checkered vertices is 
$\sum_\Omega (\rk\Omega-1)$.

The result is that $\Gamma_1\prod_\Delta\Gamma_2$ has vertices 
with black checkers, vertices with red checkers, vertices with red checkers sitting
on top of black checkers, and vertices that are completely uncheckered. Thus,
$$
\sum_\Omega (\rk\Omega-1)\leq \prod(\rk\Lambda_j-1)+N,
$$
where $N$ is the number of vertices of $\Gamma_1\prod_\Delta\Gamma_2$ that 
have a red checker but no black checker.

It remains then to estimate the number of these ``isolated'' red checkers. Observe
that a vertex of $\Gamma_1\prod_\Delta\Gamma_2$ has no black checker precisely
when it lies in the fiber, via at least one of the $\tau_j$, of a checkerless vertex
in $\Gamma_j$. Turning it around, we investigate the fibers of the checkerless 
vertices of both $\Gamma_j$. Indeed, in an $x_1$-pair, 
$$
\begin{pspicture}(0,0)(14,2)
%\showgrid
\rput(3,0){
\rput(4,1){\BoxedEPSF{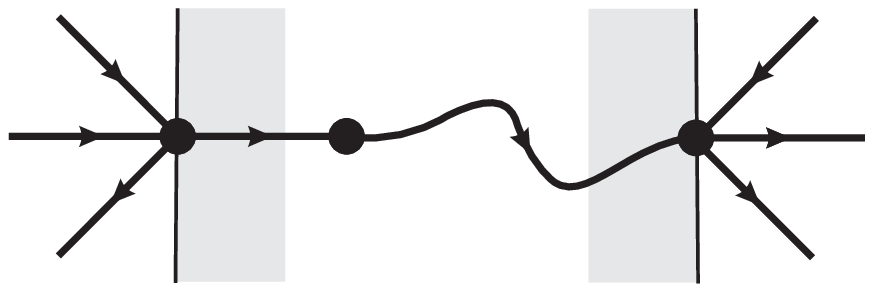 scaled 650}}
\rput(4,1.7){$\Gamma_j$}
%\rput(3.8,.9){$\gamma$}
\rput(2.8,.85){$e$}%\rput(5.8,1.25){$e_2$}
\rput(2.55,1.25){$v_1$}\rput(5.4,1.25){$v_2$}
\rput(3.6,1.25){$u$}
%\rput(1,1){$T_{e_1}$}\rput(7,1){$T_{e_2}$}
}
%\rput(14,1){$(\ddag)$}
\end{pspicture}
$$
the vertices $v_1,v_2$ and $u$ are checkerless, while $v_1,v_2$ are also checkerless
in an $x_2$-pair. We claim that no vertex in the fiber, via $t_j$, of these
five has a red checker. A vertex of $\Upsilon$ in the fiber of the boundary
vertices $v_1,v_2$ is itself a boundary vertex, hence contains no red checker.
If $v_1\times v_2\in\Upsilon$ is in the fiber of $u$ and is a boundary vertex of 
$\Upsilon$ then it carries no red checker either. If instead 
$v_1\times v_2$ is an interior vertex
then the lift to $v_1\times v_2$ of $e^{-1}$ cannot be in the coboundary
$\delta\Upsilon$, hence the terminal vertex of this lift is $\in\Upsilon$ also
and covers $v_1$. Thus, this terminal vertex is a boundary vertex for an $x_1$-pair
of $\Upsilon$, and $v_1\times v_2$ is the interior vertex from which a red
checker is removed for this pair.

The only remaining checkerless vertices of the $\Gamma_j$ unaccounted for are
those interior vertices chosen for each $x_2$-pair, and thus $N\leq$ the number
of vertices of $\Gamma_1\prod_\Delta\Gamma_2$ contained in the fibers of these.
If $u\in\Gamma_1$ is one of these interior vertices, then $u\times V_{\Gamma_2}$
are the vertices of $\Gamma_1\prod_\Delta\Gamma_2$ in the fiber. As the boundary
vertices in this fiber do not have red checkers we need only consider
the $u\times\{\text{interior vertices of }\Gamma_2\}$ with these interior vertices
precisely those of the spine $\widehat{\Lambda}_{u_2}$.
Thus our fiber is $u\times\{\text{vertices of }\widehat{\Lambda}_{u_2}\}$,
of which there are $\HH_2$, and a total of $n_{12}\HH_2$ vertices of 
$\Gamma_1\prod_\Delta\Gamma_2$ arising this way. There are also $n_{22}\HH_1$ 
vertices arising in this way from $u\in\Gamma_2$, 
and $n_{12}n_{22}$ vertices counted twice.

Thus $N\leq n_{12}\HH_2+n_{22}\HH_1-n_{12}n_{22}$, 
hence the result for $i=2$. Interchanging the checkering scheme for the $x_i$-pairs
gives the result for $i=1$.
\qed
\end{proof}

\section{Free groups and the topological dictionary}\label{free}

A group $F$ is {\em free of rank $\rk F$\/} if and only if it is 
isomorphic to the fundamental group of a connected graph of rank $\rk F$.
If\/ $\Gamma_1,\Gamma_2$ are connected graphs with 
$\pi_1(\Gamma_1,v_1)\cong\pi_1(\Gamma_2,v_2)$, then
we have $H_1(\Gamma_1)\cong H_1(\Gamma_2)$ and thus 
$\rk\Gamma_1=\rk\Gamma_2$.

The free groups so defined are of course the standard free groups and the 
rank is the usual rank
of a free group. At this stage we appeal to the existing (algebraic) theory 
of free groups, and in particular, 
that by applying Nielsen transformations, a set of generators for a free 
group 
can be transformed into a set of 
free generators whose cardinality is no greater. Thus, a finitely generated 
free group has finite rank (the converse
being obvious). From now on we use the (topologically more tractible) notion
of finite rank as a synonym for finitely generated.

Let $F$ be a free group with representation
$\varphi:F\rightarrow\pi_1(\Delta,v)$ for 
$\Delta$ connected. 
The topological dictionary is the loose term used to describe the
correspondence between algebraic properties of $F$ and topological
properties of $\Delta$ as described in 
\S\S \ref{section:topological}-\ref{topology:galoistheory}.
The non-abelian $F$ correspond to the
$\Delta$ with $\rk\Delta>1$.
A subgroup $A\subset F$ corresponds to a 
covering $p:\Lambda_u\rightarrow\Delta_v$
with $p^*\pi_1(\Lambda,u)=\varphi(A)$, and hence $\rk A=\rk\Lambda$. Thus
finitely generated subgroups correspond to finite rank $\Lambda$ and normal subgroups
to Galois coverings. Inclusion relations between subgroups correspond to 
covering relations, indices of subgroups to degrees of coverings, trivial
subgroups to simply connected coverings, conjugation to change of basepoint.

Applying the topological dictionary to the italicised results below we recover some
classical facts (see also \cites{Servatius83,Stallings83}).

\begin{enumerate}
\item \cite{Greenberg60,Karrass69}: If a finitely generated subgroup $A$ of a 
non-abelian free group $F$ 
is contained in no subgroup of infinite rank, then $A$ has finite index in $F$;
{\em Proposition \ref{finiterank:result600}}.
\item \cites{Greenberg60}: 
If a finitely generated subgroup $A$ of a free group 
$F$ contains a non-trivial normal subgroup of $F$, then it has finite index in $F$;
{\em Proposition \ref{finiterank:result700}}.
\item \cites{Burns69,Hall49}: Let $F$ be a free group, $X$ a finite subset of $F$, 
and $A$ a finitely
generated subgroup of $F$ disjoint from $X$. Then $A$ is a free
factor of a group $G$, of finite index in $F$ and disjoint from $X$;
{\em Proposition \ref{finiterank:result800}} (and the comments following it).
\item \cite{Howsen54}: If $A_1,A_2$ are finitely generated subgroups 
of a free group $F$, then the intersection of conjugates
$A_1^{g_1}\cap A_2^{g_2}$ is finitely generated for any
$g_1,g_2\in F$;
{\em Corollary \ref{pullbacks:result200}}.
\end{enumerate}

If $\Delta$ is a graph, $\rk\Delta=2$, and $A\subset F=\pi_1(\Delta,v)$, then we define
$\HH(F,A):=\HH(\Lambda_u\rightarrow\Delta_v)$ and 
$n_i(F,A):= n_i(\Lambda_u\rightarrow\Delta_v)$, where
$p:\Lambda_u\rightarrow\Delta_v$ is the covering with $p^*\pi_1(\Lambda,u)=A$. 
For an arbitrary free group $F$ with representation $\varphi:F\rightarrow\pi_1(\Delta,v)$,
define $\HH^\varphi(F,A)$ and $n^\varphi_i(F,A)$ to be
$\HH(\varphi(F),\varphi(A))$ and $n_i(\varphi(F),\varphi(A))$.

The appearance of $\varphi$ in the notation is meant to indicate that these 
quantities, unlike rank, are representation dependent. This can be both a strength
and a weakness. A weakness because it seems desirable for algebraic statements to
involve only algebraic invariants, and a strength if we have the freedom to choose
the representation, especially if the most interesting results are obtained when this
representation is not the ``obvious'' one.

For example, if $F$ is a free group with free generators $x$ and $y$, and $\Delta$ is
single vertexed with two edge loops whose homotopy classes are $a$ and $b$, then
the subgroup $A=\langle xy\rangle\subset F$ corresponds to the $\Lambda$ below
left under the obvious representation $\varphi_1(x)=a,\varphi_1(y)=b$, and
to the righthand graph via $\varphi_2(x)=a,\varphi_2(y)=a^{-1}b$:
%\begin{figure}
$$
\begin{pspicture}(0,0)(12,3)
\rput(9.5,1.5){\BoxedEPSF{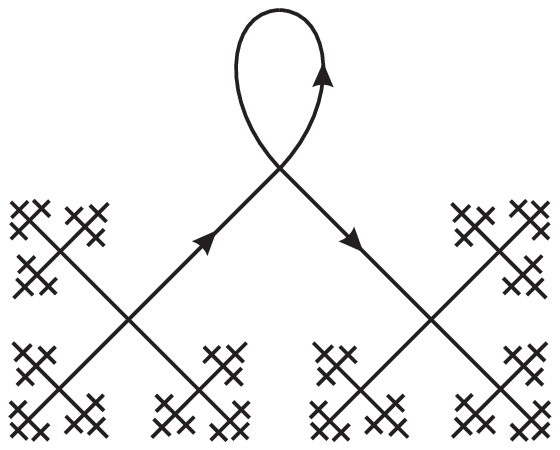 scaled 500}}
\rput(3.5,1.5){\BoxedEPSF{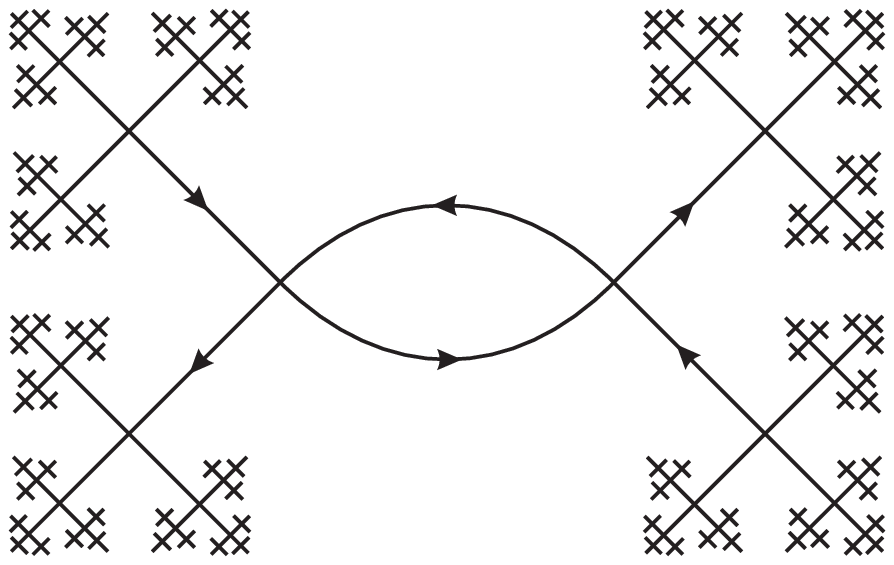 scaled 500}}
%\rput(3,1.5){\BoxedEPSF{fig18a.eps scaled 400}}
%\showgrid
\end{pspicture}
$$
Thus, $\HH^{\varphi_1}(F,A)=2,n^{\varphi_1}_{i}(F,A)=1, (i=1,2)$, 
whereas $\HH^{\varphi_2}(F,A)=1,n^{\varphi_2}_{1}(F,A)=1,n^{\varphi_2}_{2}(F,A)=0$.
%\caption{}
%\end{figure}

We now apply the toplogical dictionary to Theorem \ref{pullback:rankestimate}.
Let $\varphi:F\rightarrow\pi_1(\Delta,v)$, $A_j\subset F, (j=1,2)$, finitely generated
non-trivial subgroups, and 
$p_j:\Lambda_{u_j}\rightarrow\Delta_v, (j=1,2)$ coverings with
$\varphi(A_j)=p_j^*\pi_1(\Lambda,u_j)$.
Each non simply-connected component $\Omega$ of the pullback corresponds to some
non-trivial intersection of conjugates $A_1^{g_1}\cap A_2^{g_2}$. 
As observed in \cite{Neumann90}, these in turn
correspond to the conjugates $A_1\cap A_2^g$ for $g$ from a set of
double coset representatives for $A_2\backslash F/ A_1$.

\begin{theorem}
\label{algebraic:shn}
Let $F$ be a free group of rank two and $A_j\subset F, (j=1,2)$, finitely generated
non-trivial subgroups. Then for any representation $\varphi:F\rightarrow\pi_1(\Delta,v)$
and $i=1,2$, 
$$
\sum_g (\rk(A_1\cap A_2^g)-1)\leq 
\prod_j(\rk A_j-1)
+\HH_1\HH_2-(\HH_1-n_{1i})(\HH_2-n_{2i}),
$$
the sum over all double coset representatives $g$ for $A_2\backslash F/ A_1$
with $A_1\cap A_2^x$ non-trivial, and where $\HH_j=\HH^\varphi(F,A_j)$ and
$n_{ji}=n^\varphi_i(F,A_j)$.
\end{theorem}

This theorem should be viewed in the context of attempts %over the last half century
to prove the so-called {\em strengthened Hanna Neumann conjecture\/}: namely,
if $A_j, (j=1,2)$ are finitely
generated, non-trivial, subgroups of an arbitrary free group $F$, then 
$$
\sum_g (\rk(A_1\cap A_2^g)-1)\leq 
\prod_j(\rk A_j-1)+\varepsilon,
$$
the sum over all double coset representatives $g$ for $A_2\backslash F/ A_1$
with $A_1\cap A_2^x$ non-trivial,
where the conjecture is that $\varepsilon$ is zero, while in the existing
results, it is an error term having a long history. We 
provide a very partial, and chronological, summary of these estimates
for $\varepsilon$ in the table:
$$
\begin{tabular}{ll}
\hline\\
%Topology of graphs&Free groups\\\hline
$(\rk A_1-1)(\rk A_2-1)$&H. Neumann \cite{Neumann56}\\
$\max\{(\rk A_1-2)(\rk A_2-1),(\rk A_1-1)(\rk A_2-2)\}$&Burns \cite{Burns69}\\
$\max\{(\rk A_1-2)(\rk A_2-2)-1,0\}$&Tardos \cite{Tardos96}\\
$(\rk A_1-3)(\rk A_2-3)$&Dicks-Formanek \cite{Dicks01}\\
\\\hline
\end{tabular}
%the sum of ranks on
%the left hand side of the Theorem is $\leq \prod (\rk A_j-1)$. 
$$
(the original, unstrengthened, conjecture \cite{Neumann56} involved just the
intersection of the two subgroups, rather than their conjugates, and
the first two expressions for $\varepsilon$ were proved in this restricted
sense; the strengthened version was formulated in \cite{Neumann90}, and the 
H. Neumann and Burns estimates for $\varepsilon$ were improved to the 
strengthened case there).
Observe that as the join
$\langle A_1,A_2\rangle$ of two finitely generated subgroups 
is finitely generated, and every finitely generated free
group can be embedded as a subgroup of the free group of rank two, we may replace
the ambient free group in the conjecture with the free group of rank two.

It is hard to make a precise comparison between the $\varepsilon$ provided by 
Theorem \ref{algebraic:shn}
and those in the table. Observe that if $A_j\subset F$, with
$F$ free of rank two, then with respect to a topological representation we have
$\rk A_j=\HH_j-(n_{j1}+n_{j2})+1$. It is straight forward to find infinite families
$A_{1k},A_{2k}\subset\pi_1(\Delta,v), (k\in\Z^{>0})$, for which the error term in
Theorem \ref{algebraic:shn} is less than those in the table above for all but finitely
many $k$, or even for which the strengthened Hanna Neumann conjecture is true
by Theorem \ref{algebraic:shn}, for instance,
$$
\begin{pspicture}(0,0)(12,2)
\rput(6,1){\BoxedEPSF{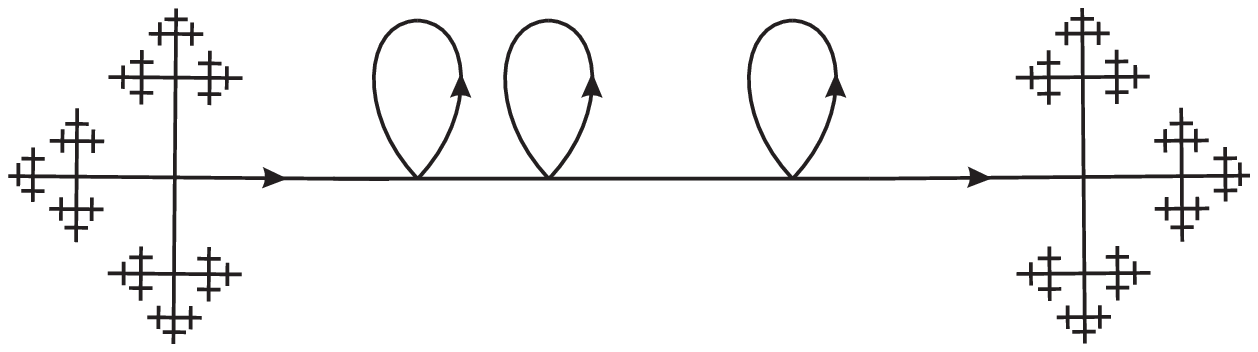 scaled 500}}
\rput{270}(5.9,.8){$\left.\begin{array}{c}
\vrule width 0 mm height 22 mm depth 0 pt\end{array}\right\}$}
\rput(5.9,.4){$k$ edge loops}
\rput(1.8,1.5){$A_{1k}=A_{2k}=$}\rput(6.2,1.5){$\ldots$}
%\showgrid
\end{pspicture}
$$
but where the error terms in the table are quadratic in $k$.

%%%%%%%%%%%%%%%%%%%%%%%%%%%%%%%%%%%%%%%%%%%%%%%%%%%%%%%%%%%%%%%%%%%%%%%%%

%\bibliography{brent}
%\bibliographystyle{plain}

%\end{bibdiv}

\end{document}